\documentclass[a4paper,english,11pt]{article}
 \usepackage[applemac]{inputenc}
\usepackage{babel}
\usepackage{epsf}
\usepackage[dvips]{epsfig}
\usepackage{subfig}
 \usepackage{graphicx}
\usepackage[colorlinks=True]{hyperref}
\usepackage{amssymb} 
\usepackage{amsmath}
\usepackage[nice]{nicefrac}
\usepackage{amsfonts}
\usepackage{dsfont}
\usepackage[T1]{fontenc}
\usepackage{wasysym}
\usepackage{amssymb}
\usepackage{stmaryrd}
\usepackage{aeguill}
\usepackage{mathrsfs}
\newtheorem{theoreme}{Th\'eor\`eme}[section]
\newtheorem{theoreme-anglais}{Theorem}[section]

\newtheorem{cor-anglais}[theoreme]{Corollary}

\newtheorem{lem}[theoreme]{Lemma}

\renewcommand{\epsilon}{\varepsilon}
\def\btab{\begin{eqnarray*}}
\def\etab{\end{eqnarray*}}
\def\beq{\begin{equation}}
\def\eeq{\end{equation}}

\newcommand{\Se}{\mathscr{S}}
\newcommand{\He}{\mathcal{H}}
\newcommand{\Kt}{M}

\newcommand{\If}{{\bf I}_{\Phi(n)}}

  \newcommand{\Sl}{\mathcal{C}}
  \newcommand{\Dl}{\mathcal{D}}
   \newcommand{\Cs}{\mathscr{C}}
\newcommand{\Ks}{\mathcal{R}}   

 \newcommand{ \un }{\mathds{1}}
 \newcommand{ \p }{\mathbb{P} }

 \newcommand{ \E }{\mathbb{E}}

\newcommand{ \Aa }{\mathscr{B}}

 \newcommand{ \R }{ \mathbb{R} }
 \newcommand{ \Z }{ \mathbb{Z} }
 \newcommand{\N}{ \mathbb{N} }

\newcommand{ \px }{\overset{\leftarrow }{x}}

\newcommand{ \pe}{ \p^{\mathcal{E}} }
\newcommand{ \Ee}{ \mathbb{E}^{\mathcal{E}} }
\newcommand{\Gg}{g}

\newcommand{ \D}{ \mathcal{D} }
\newcommand{ \A}{ \mathcal{A} }
 \newcommand{ \B}{ \mathcal{B} }

 \newcommand{ \lo }{ \mathcal{L} }

\newtheorem{The}{{\bf Theorem}}[section]
 \newtheorem{Lem}[The]{Lemma}
 
 \newtheorem{Pro}[The]{\bf Proposition}
\newtheorem{Rem}[The]{{\bf Remark}}

 \newenvironment{Pre}{\noindent \textbf{Proof.} \\ }{\hfill $
 \blacksquare$}

\makeatletter

\@addtoreset{equation}{section} \makeatother

\title{Spread of visited sites of a random walk along the generations of a branching process
\author{P. Andreoletti, P. Debs \footnote{Laboratoire MAPMO - C.N.R.S. UMR 7349 - F\'ed\'eration Denis-Poisson, Universit\'e d'Orl\'eans
(France). \newline \vspace{0.1cm}  $\quad$  MSC 2000  60J55 ;  60J80 ; 60G50 ; 60K37. \newline \vspace{0.5cm} \textit{Key words :  random walks, random environment, trees, branching random walk} }
}}

\begin{document}

\maketitle

\begin{abstract}
In this paper we consider a null recurrent random walk in random environment on a super-critical Galton-Watson tree. 
We consider the case where the log-Laplace transform $\psi$ of the branching process satisfies $\psi(1)=\psi'(1)=0$ for which G. Faraud, Y. Hu and Z. Shi in  \cite{HuShi10b}  show that, with probability one,  the largest generation visited by the walk, until the instant $n$, is of the order of $(\log n)^3$. In  \cite{AndreolettiDebs1} we prove that the largest generation entirely visited behaves almost surely like $\log n$ up to a constant. 
Here we study how the walk visits the generations $\ell=(\log n)^{1+ \zeta}$, with $0 < \zeta <2$. We obtain results in probability giving the asymptotic logarithmic behavior of the number of visited sites at a given generation. We prove that there is a phase transition at generation $(\log n)^2$ for the mean of visited sites until $n$ returns to the root. Also we show that the visited sites spread all over the tree until generation $\ell$. 

 \end{abstract}


\section{Introduction}
We start giving an iterative construction of the environment. 
Let  $(A_i,i\geq 1)$ a positive random sequence and $N$ an independent $\mathbb N$-valued random variable following a distribution $q$, in other words $\p(N=i)=q_i$ for $i\in\mathbb N$. Let $\phi$ the root of the tree and $(A(\phi^i),i\leq N_\phi))$ an independent copy of  $(A_i,i\leq N)$. Then, we draw $N_\phi$ children to $\phi$: these individuals are the first generation. Each child $\phi_i$ is associated with the corresponding $A(\phi^i)$ and so on. At the $n$-th generation, for each individual $x$ we pick $(A(x^i),i\leq N_x)$ an independent copy of $(A_i,i\leq N)$ where $N_x$ is the number of children of $x$ and $A(x^i)$ is the random variable attached to $x^i$. The set $\mathbb T$, consisting of the root and its descendants, forms a Galton-Watson tree (GW) of offspring distribution $q$ and where each vertex $x\neq\phi$ is associated with a random variable $A(x)$.\\
We denote by $|x|$ the generation of $x$, $\px$ the parent of $x$, and for convenience reasons we add $\overset{\leftarrow }{\phi}$, the parent of $\phi$. The set of environments denoted by $\textbf{E}$ is the set of all sequences  $((A(x^{i}),i\leq N_x), x \in \mathbb{T})$,  with  $P$ and $E$ respectively the associated probability measure and expectation. \\
We assume that the distribution of $(A_i,i \leq N)$ is non-degenerate and, to obtain a supercritical GW, that $E[N]>1$. Moreover we add uniform ellipticity conditions
\begin{align}
& \exists \ 0<\epsilon_0<1,\ P-a.s\  \ \forall i, \epsilon_0 \leq A_i \leq 1/ \epsilon_0, \label{hyp1}\\
&  \exists \ N_0\in\mathbb N, P-a.s\  N \leq N_0.\label{hyp2}
\end{align}
Given $\mathcal E\in \bf E$, we define a   $\mathbb T$-valued random walk $(X_n)_{n\in\mathbb N}$ starting from $\phi$ by its transition probabilities, 
\begin{align*}
& p(x,x^{i})= \frac{A(x^{i})}{\sum_{j=1}^{N_x}A(x^{j})+1},\ p(x,\px)=1-\sum_{j=1}^{N_x}p(x,x^{j}),\ p(\overset{\leftarrow }{\phi}, \phi)=1.
\end{align*}
Note that our construction implies that $(p(x,.),x\in\mathbb T)$ is an independent sequence. 
We denote by $\p^{\mathcal{E}}$ the probability measure associated to  this walk, the whole system is described under the probability $\p$, the semi-direct product of $P$ and $\p^{\mathcal{E}}$. \\
To study asymptotical behaviours associated to $(X_n)_{n\in\mathbb N}$, a quantity appears naturally: the potential process $V$ associated to the environment which is actually a branching random walk. It is defined by $V(\phi):=0$ and
$$V(x):=-\sum_{z\in\rrbracket \phi, x  \rrbracket } \log A(z),x\in\mathbb T\backslash \lbrace \phi\rbrace, $$
where $\llbracket \phi, x  \rrbracket  $ is the set of vertices on the shortest path connecting $\phi$ to $x$ and $\rrbracket \phi, x  \rrbracket=\llbracket \phi, x  \rrbracket\backslash \lbrace \phi\rbrace$.  
We put ourself in the non lattice case so $\log A_i$ can not be written as $b+c \Z$, and introduce the moment-generating function
$$\psi(t):=\log E\left[\sum_{\vert x\vert=1} e^{-tV(x)}\right],$$
characterizing the environment. Note that the hypothesis we discuss above implies that $\psi$ is defined on $\mathbb R$, and $\psi(0)>0$. In fact the hypothesis  \eqref{hyp1} and $\eqref{hyp2}$ are not always needed for our work and they could be replaced by the existence of $\psi$ in $(-\delta,1+ \delta)$ with $\delta>0$ together with the existence of a moment larger than 1 for $N$. In Section 2 for example we could lighten the hypothesis this way, but it would be much more complicated in Section 4. \\
Thanks to the work of M.V. Menshikov and D. Petritis, see \cite{MenPet}  and the first part of  \cite{Faraud} by G. Faraud, if 
\begin{align} 
\psi(1)=\psi'(1)=0 \label{hyp0}
\end{align}
then $X$ is null recurrent, with $\psi'(1)=-E\left[\sum_{\vert x\vert=1} V(x)e^{-V(x)}\right]$. In \cite{HuShi10b} (see also \cite{HuShi10a}), G. Faraud, Y. Hu and Z. Shi  study the asymptotic behavior of $\max_{ 0 \leq i \leq n}{|X_i|}=X_n^*$, i.e. the largest generation visited by the walk. Assuming \eqref{hyp0}, they prove the existence of a positive constant $a_0$ (explicitely known) such that $\p$ a.s. on the set of non-extinction of the GW
\begin{align}\label{maxx}
\lim_{n \rightarrow+\infty}\frac{X_n^*}{(\log n)^3}=a_0.
\end{align}
\noindent In \cite{AndreolettiDebs1} we were interested in the \textit{largest generation entirely visited by the walk}, that is to say the behavior of $R_n:=\sup\lbrace k\geq1,\forall  \vert z\vert= k,\lo(z,n)\geq1 \rbrace$, with $\lo$ the local time of $X$ defined by $\lo(z,n):= \sum_{k=1}^n \un_{X_k=z}$. More precisely, if \eqref{hyp0} is realized, $\p$ a.s. on the set of non-extinction  
\begin{align}\label{toutvisite}
\lim_{n\rightarrow+\infty}\frac{ R_n}{\log n}=\frac{1}{\tilde \gamma},
\end{align}
where $\tilde{\gamma}:=\sup\{a \in \R,\ \tilde J(a)>0\}$ with $\tilde{J}(a):=\inf_{t\geq 0}\{\psi(-t)-at\}$.\\
Although in \cite{AndreolettiDebs1} all recurrent cases are treated, here we focus only on the hypothesis $\eqref{hyp0}$.\\
According to \eqref{maxx} and \eqref{toutvisite}, until generation $\nicefrac{\log n}{\tilde \gamma}$  all the points are visited but $X$ does not visit generations further than $a_0(\log n)^3$. The aim of this paper is to study the asymptotic of the  number of visited sites at a given generation $(\log n)^{1+ \zeta}$ with $0< \zeta <2$. 
For this purpose we define the number of visited sites at generation $m \in \N$ until the instant $n$
$$\Kt_{n}(m):=\#\lbrace |z|=m,\mathcal L(z,n)\geq 1\rbrace,$$
and before $n$ returns to the root $K_n(.):=\Kt_{T_\phi^n}(.)$ where $T_{x}^n=\inf\{k>T_{x}^{n-1}, X_k=x\}$ for $n\geq 1$ and $T_{x}^0=0$ for $x\in\mathbb T$.\\
 \noindent Let $Z_m$ the number of descendants at generation $m \in \N$, we have $Z_1=N$. 
\noindent Our first results quantify the number of visited points at a given generation $\ell:= (\log n)^{1+ \zeta}$. 
Thanks to the hypothesis of ellipticity, $\psi$ can be written as a power series in particular, for any $x$ small enough,  $\psi(1-x)=\sum_{j=1}^{+ \infty} u_j x^j $, where $u_j=\psi^{(j)}(1)$, these are called cumulants and here $u_1:=\psi'(1)=0$, $u_2:=\psi''(1)=\sigma^2$. 
Let us define the function $f$,  for any $x$ small enough  
$$f(x):=1-\frac{x}{2 \sigma^2}+x^2 \lambda(x).$$
 $\lambda$ is the Cram\'er's  series depending on the cumulants of $\psi(1-x)$ (for more details on the Cram\'er's series see for example \cite{Petrov} p. 219-223).

\begin{The} \label{Th1}
For all $ 0  < \zeta <2$, $\epsilon>0$ independent of $\zeta$ there exists $C_0>0$ such that
\begin{align}
\lim_{n \rightarrow + \infty} \p\left(  \frac{\psi(0)}{\tilde \gamma}(1- \epsilon) \leq \frac{ \log \Kt_n(\ell)}{\log n} \leq   1-C_0\left(\frac{\log \log n}{\log n } \vee \frac{1}{(\log n)^{\zeta}}\right)  \right)=1. \label{1.6}
\end{align}
Also for all $n$ large enough, there exist two positive constants $C_1$ and $C_2$ such that
\begin{align}
  \frac{C_1 }{(\log n)^{\epsilon}}\frac{e^{(\log n) \cdot f[(\log n)^{-\zeta}]}}{(\log n)^{(1+ \tilde \zeta)/2}}  \leq {  \E[K_n(\ell)] } \leq C_2    \frac{e^{(\log n) \cdot f[(\log n)^{-\zeta}]}}{(\log n)^{(1+ \tilde \zeta)/2}}, \label{1.7}
\end{align}
with $\tilde \zeta:= \un_{0<\zeta<1}+ \zeta \un_{1 \leq \zeta<2}$.
\end{The}

\eqref{1.6} shows that, at each generation $\ell$, the cardinal of visited sites is at least $n^{\psi(0)(1-\epsilon)/ \tilde \gamma}$ for any $\zeta$, that is to say like the last generation entirely visited $R_n$ ($\psi(0)/ \tilde \gamma<1$, by convexity of $\psi$ and the fact  that $\psi(1)=0$). Also the upper bound of $ \Kt_n(\ell)$ is at most of the order of $n e^{-C_3 (\log n)^{1-\zeta}}/ (\log n)^{C_4}$, with $C_3,C_4>0$. This suggests that it may have a phase transition when $\zeta=1$. Although we are not able to show this for $\Kt_n(\ell)$ the existence of a phase transition is proved in \eqref{1.7} for  the mean of $K_n(\ell)$. 
Indeed by definition of $f$,
$$ (\log n) f[ (\log n)^{-\zeta}]=\log n-\frac{(\log n)^{1-\zeta}}{2 \sigma^2}+(\log n)^{1-2\zeta} \lambda((\log n)^{-\zeta}) $$
We can see that in the neighborhood of generation $(\log n)^2$ that is to say when $\zeta=1$, the asymptotic behavior of $\mathcal{N}_{\zeta}:=\E[K_n(\ell)]$ changes. We easily check that for all $0 <\zeta< \zeta' \leq 1$, $\lim_{n \rightarrow + \infty} \mathcal{N}_{\zeta'}/\mathcal{N}_{\zeta}=+ \infty$ 
whereas for all $1 \leq \zeta< \zeta' <2$, $\lim_{n \rightarrow + \infty} \mathcal{N}_{\zeta'}/\mathcal{N}_{\zeta}=0$. 
So the generations of order $(\log n)^2$ are, in mean, the most visited generation (in term of distinct site visited) until $n$ returns to the origin. Finally notice that when $\zeta>1/2$ we are in a Gaussian behavior as $ e^{(\log n) f[ (\log n)^{-\zeta}]} \sim ne^{-\frac{(\log n)^{1-\zeta}}{2 \sigma^2}}$, and when $\zeta \geq 1$, $ e^{(\log n) f[ (\log n)^{-\zeta}]}  \sim  n $.

In order to establish our second result, recall Neveu's notation to introduce a partial order on our tree. In \cite{Neveu}, to each vertex $x$ at generation $m\in\mathbb N$, Neveu associates a sequence $x_1\dots x_m $ where $x_i\in \mathbb N$, to simplify we write $x=x_1\dots x_m $.\\
This sequence gives the complete ``genealogy'' of $x$:  if $y=x_1\dots x_i$ with $|y|=i<m$, $y$ is the unique ancestor of $x$ at generation $i$ and we write $y<x$.\\
For instance $\overset{\leftarrow}{x}=x_1\dots x_{m-1}$ and $1\leq x_m\leq N_{\overset{\leftarrow}{x}}$, in other words $x$ is the $x_m$-th child of $\overset{\leftarrow}{x}$.\\
To extend this partial order for $|x|=|z|$, we write $x<z$ if there exists $i<m$ such that $x_k=z_k$ for $k<i$ and $x_i<z_i$. Hence we can number individuals at a given generation ``from the left to the right'' and for $A$ a subset of $\lbrace z\in\mathbb T, |z|=m\rbrace$,  ${\bf inf}\,A$ and ${\bf sup}\, A$ are respectively the minimum and maximum associated to this numbering.\\
Our last result gives an idea of the way the visited points spread on the tree, for this purpose we introduce clusters:
let $z \in \mathbb T$ and $m \geq \vert z\vert$, we call cluster  issued from $z$ at generation $m$ denoted $ \mathcal{C}_m(z)$, the set of descendants $u$ of $z$ such that  $\vert u\vert=m $, in other words
\begin{align}
\mathcal{C}_m(z):= \{u>z,\ |u|=m \}. \label{defclust}  \end{align}
At some point we need to quantify the number of individuals between two disjoint clusters with common generations. For given initial and terminal generations, denote $\Cs$ a set of disjoint clusters. 
Let $(\Dl_j, 1\leq j \leq |\Cs|)$, with $|\Cs|$ the cardinal of $\Cs$,  an ordered  sequence of (disjoint) clusters belonging to $\Cs$, that is to say for all $j$, $ {\bf sup} \, \Dl_j   < {\bf inf} \,\Dl_{j+1}  $. We define the minimal distance between clusters in the following way $ {\bf{D}}(\Cs):= \min_{ 1 \leq j \leq |\Cs|-2 }({\bf inf} \, \Dl_{j+2}  -{\bf sup} \,  \Dl_{j}  ) $, where, by definition, ${\bf inf} \, \Dl_{j+2}  -{\bf sup} \   \Dl_{j} $ is the number of individuals between ${\bf sup} \, \Dl_{j}$ and ${\bf inf} \,  \Dl_{j+2}$. Notice that we do not look at two successive clusters, but  two successive separate by one. We now state a second result
\begin{The} \label{Th2}
For $ 0 < \zeta <2$ and $\epsilon>0$ recalling that $\ell= (\log n)^{1+ \zeta}$
\begin{align}
& \lim_{n \rightarrow + \infty} \p\left( \max_{|z|=\ell- \log n/\tilde{\gamma} } \min_{ y \in \mathcal{C}_{\ell} (z)} {\lo(y,n) \geq 1}   \right) =1, \label{1.8} \\
& \lim_{n \rightarrow + \infty} \p\left(\min_{z \in \mathcal{C}_{\epsilon \ell ^{1/3}} (\phi)}  \max_{y>z,|y|=\ell } \lo(y,n) \geq 1  \right) =1. \label{1.9}
\end{align}
Let $k_n$, $h_n$ and $r_n$ positive sequences of integers such that $k_n r_n+(k_n-1) h_n=\ell$. For all $ 1 \leq i \leq k_n $, let us denote $\Cs_i$ a set of clusters initiated at generation $ (i-1) (r_n+ h_n) $  and with end points at generation $ i r_n+(i-1) h_n $ (see Figure \ref{fig3}), also define the following event for all $m>0$ and $q>0$
\begin{align*}
\mathscr{A}_i(m,q):= \bigcup_{\Cs_i}\left\{ \left\lbrace{ |\Cs_i| \geq q,\bf{D}}(\Cs_i) \geq m\right\rbrace \bigcap_{ \Dl \in \Cs_i } \left\{ {\forall z\in \Dl, \lo(z,n)\geq 1} \right\}   \right\} .
\end{align*}
 There exist $0<{\bf k }<1 \wedge \zeta$, $0<{\bf r}<1$ with $ 0<{\bf k}+{\bf r} \leq 1$  and for $k_n=(\log n)^{{\bf k}}$, $r_n=(\log n)^{{\bf r}}$
\begin{align}
& \lim_{n \rightarrow + \infty} \p\left(\bigcap_{i=2}^{k_n}\mathscr{A}_i(e^{ \psi(0)h_n/2},e^{\psi(0)r_n(i-1)/2}) \right)=1 \label{1.10b}.
\end{align}

\end{The}

\begin{figure}[h!]
\begin{center}
\includegraphics[width=5cm]{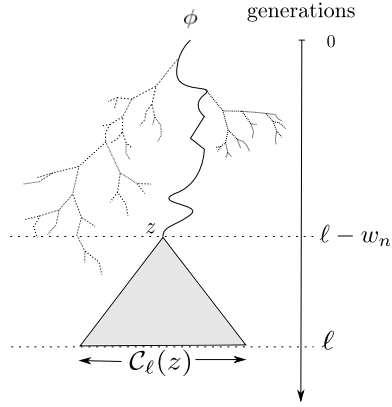}
\end{center}
\caption{One large cluster} \label{fig1}
\end{figure}

\noindent \eqref{1.8} implies the existence of a cluster starting at a generation $\ell- \log n/ \tilde \gamma $ completely visited (see Figure \ref{fig1}). As conditionnaly on the tree until generation $|z|$, $|\mathcal{C}_{\ell} (z)|$ is equal in law to $Z_{\ell-|z|}=Z_{\psi(0) \log n/\tilde{\gamma}} $, this cluster is large and, in particular,  \eqref{1.8}  implies the lower bound in \eqref{1.6}. \\
\eqref{1.9} tells that we can find visited individuals at generation $\ell= (\log n)^{1+ \zeta}$, with a common ancestor to a generation close to the root, that is to say before generation $\epsilon \ell ^{1/3}$ (see Figure \ref{fig4}). Thus, with a probability close to one, at least $e^{\epsilon (1- \epsilon) \psi(0) \ell^{1/3}/2 }$ individuals of generation $\ell$ separate by at least $e^{\psi(0) \ell /2}$ individuals   of the same generation $\ell$, are visited. \\
Finally \eqref{1.10b} tells that if we make cuts regularly on the tree we can find many visited clusters (which number increases with the generation) well separated. In particular these visited clusters can not be in a same large visited clusters as they are separated by at least $e^{\psi(0) h_n/2} \sim e^{\psi(0)(\log n)^{1+ \zeta-\bf{k}}/2}>n$ individuals (see also Figure \ref{fig3}). 

\noindent \\ To obtain these results we show that $K_{n}(\ell)$ can be linked to a  random variable depending only on the random environment and $n$. For all $z\in\mathbb T$, all integer $k$ and all real $a$, we define the random variable 
$$\Ks_{a}^z(k):= \#\lbrace u>z, |u|=k, \overline{V}(z) \leq a\rbrace,$$
where $\overline {V}(z)=\max_{u \in \rrbracket \phi, z \rrbracket } V(u)$. For notational simplicity, we write $\Ks_{a}(k)$ for $\Ks_{a}^\phi(k)$.  
 We obtain the following

\begin{Pro} \label{ExpKr} Let $\epsilon>0$ and $\Phi$ a sequence such that
\begin{align}
 (1- \epsilon)\log n \leq \Phi(n) \leq \log n+o(\log n). \label{PHI} \end{align}
Then, for all $0 < \zeta <2$ there exists $C_0'>0$
\begin{align}
& \lim_{n \rightarrow + \infty} P\left(  \frac{\psi(0)}{\tilde \gamma}(1- \epsilon) \leq \frac{ \log \Ks_{\Phi(n)}(\ell) }{\Phi(n)} \leq   1-C_0'\left(\frac{\log \log n}{\Phi(n)} \vee \frac{\Phi(n)}{\ell }\right)     \right)=1, \label{Prop2} \\
& E[\Ks_{\Phi(n)}(\ell)] \asymp(\ell^{-1} \un_{0<\zeta<1}+\Phi(n) {\ell^{-3/2}}\un_{1\leq \zeta<2}){e^{\Phi(n)f(\Phi(n)/ \ell)}}. \label{Prop1}
\end{align}
\end{Pro}

We use the notation $a_n \asymp b_n$ when there exists two positive constants $c_1$ and $c_2$ such that $ c_1 b_n \leq  a_n \leq c_2 b_n$ for all $n$ large enough.
The lack of precision for the first result shows no difference between 
$\Ks_{\log n}(\ell)$ and  $\Kt_{n}(\ell)$ (see \eqref{1.6}), unlike between the means of  $\Ks_{\log n}(\ell)$ and ${K}_{n}(\ell)$.

\noindent \\ The rest of the paper is organized as follow: in Section 2 we study $\Ks_{\Phi(n)}(\ell)$ and prove Proposition \ref{ExpKr}. In Section 3 we link $\Ks_{\Phi(n)}(\ell)$ and $\Kt_{n}(\ell)$, which leads to Theorem \ref{Th1} and \eqref{1.8} of Theorem \ref{Th2}. In Section 4  we prove the end of Theorem \ref{Th2}. Also we add an appendix where we state known results on branching processes and local limit theorems for sums of i.i.d. random variables.\\
Note that for typographical simplicity, we do not distinguish a real number and its integer part throughout the article.

\section{Expectation and bounds of  $\boldsymbol{\Ks_{\Phi(n)}(\ell)}$}

In this section we only work with the environment more especially with what we call \textit{number of  accessible points} $\Ks_{\Phi(n)}(\ell)$.
\subsection{Expectation of  $\boldsymbol{\Ks_{\Phi(n)}(\ell)}$ (proof of (\ref{Prop1}))}

\noindent According to Biggins-Kyprianou identity (also called many-to-one formula, see part A of appendix), $ E[\Ks_{\Phi(n)}(\ell)]=E \left[e^{S_\ell} \un_{\bar{S}_\ell  \leq \Phi( n)}\right]$ where $S_j$ is a centered random walk,
we only have to prove

\begin{Lem} \label{lem2.1} For all $\epsilon>0$, $\Phi$ satisfying  \eqref{PHI}, for all $0 < \zeta<2$
$$
E \left[e^{S_\ell} \un_{\bar{S}_\ell  \leq \Phi( n)}\right]  \asymp e^{\Phi(n)f(\Phi(n)/\ell)} \left[\ell^{-1}\un_{0<\zeta \leq 1} + \Phi(n) \ell^{-3/2 } \un_{1\leq \zeta<2}\right]. $$
\end{Lem}

\begin{Pre}
For $\epsilon>0$: 
$$ E \left[e^{S_\ell} \un_{ \bar{S}_\ell \in \If }\right] \leq E \left[e^{S_\ell} \un_{\bar{S}_\ell  \leq \Phi( n)}\right] \leq E \left[e^{S_\ell} \un_{ \bar{S}_\ell  \in \If}\right]+ e^{\Phi(n)(1- \epsilon)}.$$
with $I_{\Phi}:=]\Phi( n) (1- \epsilon),\Phi( n)]$. For  every sequence $(u_n)_{n\in \N}$, we denote $\bar u_j:= \max_{ 1 \leq i \leq j } u_i$ and $\underline u_j := \min_{1 \leq i \leq j} u_i $, also let $\Se_j:=\lbrace \bar S_{j-1}<S_j=\bar S_\ell\rbrace$.  First, as $S_0=0$
$$\sum_{j=1}^\ell E \left[e^{S_\ell} \un_{{S}_j  \in \If,\Se_j}\right] \leq E \left[e^{S_\ell} \un_{\bar{S}_\ell  \in \If}\right] \leq \sum_{j=1}^\ell E \left[e^{S_\ell} \un_{{S}_j  \in \If,\Se_j}\right]+1.$$
For $0\leq i\leq j$, let $\tilde S_i:=S_j-S_{j-i}$, with this notation $ \{\bar S_{j-1}<S_{j}  \}=\{ \tilde{\underline{S}}_{j-1} >0 \}$ and $\tilde S_j=S_j$. Writing $S$ as a sum of i.i.d. random variables, we easily see that $(S_i)_{0\leq i\leq j}$ and $(\tilde S_i)_{0\leq i\leq j}$ have the same law. Then, conditioning on $\sigma\{S_k,k\leq j \}$ 
 \begin{align}
 E \left[e^{S_\ell} \un_{{S}_j  \in \If,\Se_j}\right]= D_j  F_{\ell-j}   \label{eq1.0}
 \end{align}
with $F_{m}:=E\left[e^{S_{m}}\un_{\bar S_{m \leq 0}} \right]$  and $D_j:=E\left[e^{S_j} \un_{S_j\in \If, \underline{S}_{j-1}>0} \right]$. \\
By \eqref{F}, $\forall j\leq \ell,\, F_{\ell-j}\asymp  (\ell-j+1)^{-3/2}$ then it remains to estimate $D_j$. For any $A>0$ 
 \begin{align*}
 D_j  & = \sum_{k=\Phi(n)(1- \epsilon)+1}^{ \Phi(n)} E\left[e^{S_j} \un_{k-1 <  S_j \leq  k, \underline{S}_j>0} \right] \\ &\asymp  
  \sum_{k=\Phi(n)(1- \epsilon)+1}^{\Phi(n)}e^k P( k-1 <  S_j \leq  k, \underline{S}_j>0) (\un_{k \leq A j^{1/2}}+\un_{k > A j^{1/2}} ) \\
&=:D_j^1+D^2_j.
%
 \end{align*} 
We now need to distinguish the cases $0 < \zeta < 1 $ and $1 \leq \zeta < 2 $. \\
\textit{When $0 < \zeta < 1 $}, $D^1_j=0$ as $k > A\ell^{1/2}$. Also using Lemma \ref{beautiful} for $A$ large enough
\begin{align}
H^2_j& := \sum_{k=\Phi(n)(1- \epsilon)+1}^{\Phi(n)}e^k P( k-1 <  S_j \leq  k, \underline{S}_j>0) \un_{ A j^{1/2}<k< \epsilon j}  \nonumber \\
 & \asymp \sum_{k=\Phi(n)(1- \epsilon)+1}^{\epsilon j \wedge \Phi(n)}e^k P( k-1 <  S_j  \leq   k, \underline{S}_j>0)   \asymp \sum_{k=\Phi(n)(1- \epsilon)+1}^{\epsilon j \wedge \Phi(n)}  \frac{e^{k f(k/j)}}{j}, \label{2.3bm}
\end{align}
recall that $f(x)=1-x/(2 \sigma^2 )+x^2\lambda(x)$ where $\lambda$ is the Cram\'er's serie associated to $V$. 
\eqref{2.3bm} implies that $H^2_j \geq c_-e^{(\epsilon j \wedge \Phi(n)) f( (\epsilon j \wedge \Phi(n))/j)} / j $. For the upper bound, we can assume without loss of generality that $\varepsilon$ is small enough to ensure that for $|x|\leq \varepsilon$, $\lambda(x)$ converges and $f^\prime (x)$ is negative.  Therefore, the derivative of $F$ defined by $F(x):=e^{xf(x/j)}/f(x/j)$ satisfies in the same interval  $F'(x)\geq e^{xf(x/j)}  - c_+ x  e^{xf(x/j)}/j \geq e^{xf(x/j)}(1  - c_+ \epsilon)    $, with $c_+>0$.    Integrating this last inequality, for $\epsilon$ small enough
\begin{align}
  \int_{\Phi(n)(1- \epsilon)}^{\epsilon j \wedge \Phi(n)} e^{xf(x/j)} \leq \left[F(x) \right]_{\Phi(n)(1- \epsilon)}^{\epsilon j \wedge \Phi(n)}/(1-c_+ \epsilon)  
  \leq C_+ e^{(\epsilon j \wedge \Phi(n))f( (\epsilon j \wedge \Phi(n))/j)},
\end{align}
and finally 
\begin{align*}
H^2_j  & \leq \frac{C_+}{j}  \int_{\Phi(n)(1-\epsilon)}^{\Phi(n)} e^{xf(x/j)}  dx 
 \leq \frac{C_+}{j} e^{(\epsilon j \wedge \Phi(n))f((\epsilon j \wedge \Phi(n))/j)}.
\end{align*}
Note that for $s>0$ small enough,  $\psi(1-s) \leq  {s}/{2}$. So the exponential Markov inequality applied to  $P(sS_j>sk)$ and the identity $E[e^{s S_j}]=e^{j\psi(1-s)}$ yield
\begin{align*} 
\tilde H^2_j:=\sum_{k=\Phi(n)(1- \epsilon)+1}^{ \Phi(n)} e^k P( k-1 <  S_j \leq   k, \underline{S}_j>0) \un_{ k \geq \epsilon j   } \leq e^{\Phi(n)(1-  s/2)}. 
\end{align*} 
In particular $\tilde H^2_j=o(H_j^2)$ for all $j\geq \Phi(n)^{1+u}$ with $u>0$. Finally, as for any $\epsilon$ small enough $(\epsilon j \wedge \Phi(n))f( (\epsilon j \wedge \Phi(n))/j)$ is increasing in $j$, for any $0<u<\zeta$  
 \begin{align} \sum_{j=1}^{\ell} D_j F_{\ell-j} \asymp \sum_{j= \Phi(n)^{1+u}}^{\ell}  (\ell-j+1)^{-3/2} j^{-1} e^{(\epsilon j \wedge \Phi(n))f((\epsilon j \wedge \Phi(n))/j)}  \asymp e^{ \Phi(n) f( \Phi(n)/ \ell)}/ \ell . \label{2.5d} \end{align}
 Indeed, writing $\sum_{j= \Phi(n)^{1+u}}^{\ell}  (\ell-j+1)^{-3/2} j^{-1} e^{(\epsilon j \wedge \Phi(n))f((\epsilon j \wedge \Phi(n))/j)} :=\sum_{j= \Phi(n)^{1+u}}^{\ell} G_j$, $\sum_{j= \Phi(n)^{1+u}}^{\ell}G_j \geq G_\ell$ and as 
 \begin{eqnarray*}
 \sum_{j= \Phi(n)^{1+u}}^{\frac{\ell}{2}}G_j&\leq& \frac{C_+e^{ \Phi(n) f( \Phi(n)/ \ell)}}{\phi(n)^{1+u}\ell^{\frac{1}{2}}}\leq \frac{C_+ e^{ \Phi(n) f( \Phi(n)/ \ell)}}{ \ell},\\
  \sum_{j= \frac{\ell}{2}+1}^{\ell}G_j&\leq&  \frac{C_+ e^{ \Phi(n) f( \Phi(n)/ \ell)}}{ \ell} \sum_{j= \frac{\ell}{2}+1}^{\ell}{(\ell-j+1)^{-\frac{3}{2}}}\leq \frac{C_+ e^{ \Phi(n) f( \Phi(n)/ \ell)}}{ \ell},
\end{eqnarray*}
\eqref{2.5d} follows.\\
\textit{When $1 \leq \zeta < 2 $}, we prove that the main contribution comes from $D^1_j$. As for any $n$ large enough,
 $\Phi(n) \leq A \ell^{1/2}$ for some $A>0$, for any $j \geq  (\Phi(n)/A)^2$ using  Lemma \ref{beautiful}
\begin{align*}
D^1_j \asymp \sum_{k=\Phi(n)(1- \epsilon)}^{ \Phi(n)}  \frac{k e^k}{ j^{3/2}}e^{-k^2/(2 \sigma^2 j)}\asymp  \frac{\Phi(n)}{j^{{3}/{2}}} e^{\Phi(n)f(\Phi(n)/\ell)}.
\end{align*}

\noindent When $j<(\Phi(n)/A)^2$, similar computations than for $H_j$ and $\tilde H_j$ give \begin{align*} \sum_{j=1}^{\ell}D^2_jF_{\ell-j} & \leq \sum_{j=1}^{(\Phi(n)/A)^2}F_{\ell-j}(H_j^2+ \tilde H_j^2)  
 \leq C_+ e^{\Phi(n)(1-(A^2/2 \sigma^2 \wedge s/2) )}.
\end{align*}
Finally  
$\sum_{j=1}^\ell D_j F_{\ell-j} \asymp  {\Phi(n)   e^{ \Phi(n) f( \Phi(n)/ \ell)} }{\ell^{-3/2}},$
this together with \eqref{2.5d} finishes the proof.
\end{Pre}

\subsection{Bounds for $ \boldsymbol{\log \Ks_{\Phi(n)}(\ell)}$ (proof of (\ref{Prop2}))}

\noindent The upper bound is a direct consequence of  Markov inequality and \eqref{Prop1}.\\
For the lower bound, we first need an estimation on the deviation of $\min_{|z|=m}  \overline{V}(m)$, this topic has been studied in details in \cite{HuShi10b}, 

\begin{Pro}  \label{lem2.4} 
Let $a_n$ a positive sequence such that  $a_n \sim n^{1/3} $, there exists $b_0>0$ such that for any $ 0<b< b_0$
\begin{align}
\lim_{n \rightarrow + \infty}  \frac{1}{a_n} \log P\left(\min_{|z|=n} \overline{V}(z) \leq b a_n\right)=b-b_0. \label{2.11}
\end{align}
\end{Pro}

\noindent A useful consequence of the above Proposition is the following

\begin{Lem} \label{ProD3} 
Assume that $a_n$ is a positive increasing sequence such that  $a_n \sim n^{1/3} $,  there exists a constant $\mu>0$ such that for any $n$ large enough
\begin{align}
 P\left(\min_{|z|=n} \overline{V}(z) > \mu a_n\right)\leq  \lambda_n+o(\lambda_n), \label{2.9}
\end{align}
where $\lambda_n=e^{-c_1e^{c_1  a_n}}$ if $q_0+q_1= 0$, and $\lambda_n=e^{-{c_1  a_n}}$ otherwise, also $c_1>0$ depends only on the distribution $P$.  
\end{Lem}

\begin{Pre} Clearly for $z_1<z$, $\overline{V}(z) \leq \overline{V}(z_1) + \tilde V(z_1,z)$ where $ \tilde V(z_1,z)=\max_{z_1 < x \leq z}V(x)-V(z_1)$. In the sequel, writing $ \tilde V(z_1,z)$ implies that $z_1<z$ implicitly. For $ 0<\eta<1$ and $v_n:= \nicefrac{\eta b_0 a_n}{ \alpha}$ 
\begin{align*}
 P\left(\min_{|z|=n} \overline{V}(z) > 2 \eta b_0  a_n\right) \leq P\left(\min_{|z_1|=v_n}\min_{|z|=n} \tilde{V}(z_1,z) > 2\eta b_0  a_n-\max_{|z_1|=v_n}\overline{V}(z_1) \right).
\end{align*}
Using that $\max_{|z_1|=v_n}\overline{V}(z_1) \leq \alpha v_n$ by ellipticity and for $\mathcal{A}_n:= \left\{ Z_{v_n} \geq e^{ \eta \psi(0) v_n }\right\}$
\begin{align*}
 P\left(\min_{|z|=n}\overline{V}(z) > 2\eta b_0a_n\right) & \leq P\left(\min_{|z_1|=v_n}\min_{|z|=n} \tilde{V}(z_1,z) > \eta b_0{a}_n \right)  \\
 &  \leq  P\left(\min_{|z_1|=v_n}\min_{|z|=n} \tilde{V}(z_1,z) > \eta b_0 a_n ,\mathcal{A}_n \right) +P(\bar {\mathcal{A}}_n).
\end{align*}
Theorem \ref{ND} tells that if $q_0+q_1 > 0$, there exists $\nu>0$ such that  $P(\bar {\mathcal{A}}_n)\leq e^{- \nu (1- \eta) \psi(0) v_n }$, otherwise there exists $\beta'>0$ such that $\log P(\bar {\mathcal{A}}_n) \sim -e^{\beta'(1- \eta) \psi(0) v_n} $. 
Stationarity gives that $ \min_{|z|=n} \tilde{V}(z_1,z)$ and  $\min_{|z|=n-v_n} \overline{V}(z)$ have the same law, and independence of the sub-branching processes rooted at generation $v_n$  together with \eqref{2.11} imply
\begin{align*}
 P\left(\min_{|z_1|=v_n}\min_{|z|=n} \tilde{V}(z_1,z) > \eta b_0 a_n, {\mathcal{A}}_n \right) 
& \leq  P\left(\min_{|z|=n-v_n} \overline{V}(z) >\eta b_0 a_n \right)^{e^{\eta \psi(0)  v_n}} \\
& \leq \left( 1- e^{-(b_0(1- \eta)+o(1))a_n } \right)^{e^{ \eta  \psi(0)v_n}},
\end{align*}
we conclude choosing $\eta$ sufficiently close to $1$ to get $(1- \eta)< \eta^2 \psi(0)/ \alpha$.
\end{Pre}


To obtain the lower bound for $\log \Ks_{\Phi(n)}(\ell)$,  we prove the existence of a cluster $\mathcal{C}_{\ell}(z)$ (see \eqref{defclust}) with $|z|=\ell-w_n$ where $w_n:= \Phi(n)(1- \epsilon)/\tilde \gamma $ and such that $\forall z^\prime\in \mathcal{C}_{\ell}(z),\,   \overline V(z') \leq \Phi(n)$. In other words for $|z|<\ell$, let  $Z^z_\ell$ the number of descendants of $z$ at generation $\ell$, we prove
\begin{align} \lim_{n \rightarrow + \infty} P\left(\bigcup_{|z|=\ell-w_n } \left\{\# \lbrace z' \in \mathcal{C}_{\ell}(z),\,\overline{V}(z') \leq \Phi(n)\rbrace=Z^z_{\ell}\right\}\right)=1, \label{2.7lim}
\end{align}
 which implies according Theorem \ref{ND} that
$$
\lim_{n \rightarrow + \infty} P \left(\Ks_{\Phi(n)}(\ell)\geq e^{\psi(0) w_n(1- \epsilon)}\right)=1.
$$
Let $\mathcal{B}:=\bigcup_{|z|=\ell-w_n} \{\overline{V}(z) \leq  y_n,\Ks_{\Phi(n)-y_n}^z(\ell)=Z^z_{\ell} \}$  where 
$y_n:= \mu \ell^{1/3}$
\begin{align*}
& P( \mathcal{B})
 \geq P\left(\left\lbrace \Ks_{y_n}(\ell-w_n)\geq 1\right\rbrace\bigcap\bigcup_{|z|=\ell-w_n, \overline{V}(z)\leq y_n}\left\lbrace \Ks^z_{\Phi(n)-y_n}(\ell)=Z^z_{\ell} \right\rbrace  \right)\\
&=\sum_{k\geq1}P\left. \left(\Ks_{y_n}(\ell-w_n)=k\right)P\left(\bigcup_{|z|=\ell-w_n,\overline{V}(z)\leq y_n}\left\{ \Ks^z_{\Phi(n)-y_n}(\ell)=Z^z_{\ell} \right\rbrace \right| \Ks_{y_n}(\ell-w_n)=k     \right).
\end{align*}
Let us denote $z_1,\dots, z_k,\dots$ the ordered  points at generation $\ell-w_n$ satisfying $\overline{V}(z_i) \leq y_n $. Conditionally on $ \lbrace\Ks_{y_n}(\ell-w_n)=k  \rbrace$,  $z_1$ exists and 
$$\left\lbrace \Ks^{z_1}_{\Phi(n)-y_n}(\ell)=Z^{z_1}_{\ell}\right\rbrace\subset \bigcup_{|z|=\ell-w_n, \overline{V}(z)\leq y_n}\left\lbrace \Ks^z_{\Phi(n)-y_n}(\ell)=Z^{z}_{\ell}\right\rbrace .$$ 
Furthermore, by stationarity $\Ks^{z_1}_{\Phi(n)-y_n}(\ell)$ and $ \Ks_{\Phi(n)-y_n}(w_n)$ have the same law, so
\begin{align*}
P(\mathcal B) &
\geq P\left(\Ks_{\Phi(n)-y_n}(w_n) = Z_{w_n}\right) \sum_{k\geq1}P\left(\Ks_{y_n}(\ell-w_n)=k\right)\\
&\geq P\left(\max_{|z|=w_n} \overline{V}(z) \leq \Phi(n)-y_n \right)
 P\left(\min_{|z|=\ell-w_n} \overline{V}(z) \leq  y_n \right)
\end{align*}
As $y_n=o(\Phi(n))$, the first probability tends to one thanks to a result of Mac-Diarmid \cite{McDiarmid} (see also \cite{AndreolettiDebs1} Lemma 2.1), so does the second one as a consequence of Lemma  \ref{ProD3}. \hfill $\blacksquare$


\section{Expectation of $\boldsymbol{K_n(\ell)}$,  bounds for $\boldsymbol{\log K_n(\ell)}$ and $\boldsymbol{\log \Kt_n(\ell)}$}

\subsection{ Proof of (\ref{1.7}) }

We start with general upper and lower bounds for the annealed expectation of $K_n(\ell)$.



\begin{lem}\label{ultimlem} For $n\in\mathbb N$:
\begin{align*}
C_-  (n A^-_n+ B^-_n) \leq \E[K_n(\ell)] \leq C_+ (n A^+_n+B^+_n)
\end{align*}
where
\begin{align*}
A_n^+:=E \left[e^{S_\ell- \bar S_\ell} \un_{\sum_{i=1}^\ell e^{S_i} > c_- n}\right], 
B_n^+:=E \left[e^{S_\ell} \un_{\bar{S}_\ell  \leq  \log(c_+ n)}\right], \\
A^-_n:=E \left[\frac{e^{S_\ell}}{\sum_{i=1}^\ell e^{S_i} } \un_{\bar{S}_\ell > \log ( c_+n)}\right] \textrm{ and } B^-_n:= E \left[e^{S_\ell} \un_{\sum_{i=1}^\ell e^{S_i} \leq  c_- n }\right],
\end{align*}
 $C_-$ and $c_-$ (respectively $C_+$ and $c_+$) are positive constants that may  decrease (respectively increase) from line to line.
\end{lem}

\begin{Pre}
Markov property gives $\Ee[K_n(\ell)] = \sum_{\vert z\vert =\ell} (1- e^{n \log(1- p_z) })$, with $p_z:=\pe_\phi(T_z<T_\phi)$. Obviously on $\lbrace np_z\geq 1\rbrace$, $  1-e^{-1}\leq 1- e^{n \log(1-p_z) } \leq 1.$
As for $x\in[0;1[$, $-x(1+\nicefrac{x}{2})\leq \log(1-x)\leq -x$ and $x(1-\nicefrac{x}{2})\leq1-e^{-x}\leq x $, on $\lbrace np_z<1\rbrace$
$$-3p_z/2 \leq \log (1-p_z) \leq -p_z\mbox{ and }\frac{np_z}{4}\leq 1- e^{n \log(1-p_z) }\leq\frac{3np_z}{2},$$
then
\begin{align*}
C_-({n p_z} \un_{n p_z <1}+ \un_{n p_z  \geq 1})  \leq 1- e^{n \log(1-p_z) } \leq  C_+( n p_z \un_{n p_z <1}+ \un_{n p_z  \geq 1}). 
\end{align*}
Using successively the fact that  $c_-(\sum_{x\in\rrbracket \phi,z\rrbracket}e^{ V(x)})^{-1} \leq p_z \leq c_+ e^{- \overline{V}(z)}$ and  Biggins-Kyprianou identity (see Appendix \ref{apA})
\begin{eqnarray*}
B_n^-&    \leq E\left[  \sum_{\vert z\vert =\ell} \mathds{1}_{np_z \geq 1}\right]\leq &B_n^+.
\end{eqnarray*}
Similar arguments show $C_-A_n^-\leq E\left[  \sum_{\vert z\vert =\ell}p_z \mathds{1}_{np_z < 1}\right]\leq C_+ A_n^+$.\end{Pre}

\noindent We now give upper bounds for  $B_n^+$ and $A_n^+$, and a lower bound for $A_n^-$. \\ 
$\bullet$ \textit{For} $B_n^+$, we use Lemma \ref{lem2.1} taking $\Phi(n)= \log n$. \\
$\bullet$ \textit{For} $A_n^+$, first note that
$\lbrace  \sum_{i=1}^\ell e^{S_i} > c_- n\rbrace\subset\left\lbrace\bar S_\ell>d_n\right\rbrace $, with $d_n=\log({c_-n}/{\ell})$. Recalling the arguments given in \eqref{eq1.0}, $A_n^+$ is bounded from above by
\begin{align}  E \left[e^{S_\ell- \bar S_\ell} \un_{ {\bar S_\ell} > d_n}\right] 
=  \sum_{j=\nicefrac{d_n}{\alpha}}^\ell P\left({ S_j > d_n, \underline{S}_j>0}\right) E\left[e^{S_{\ell-j}} \un_{\bar{S}_{\ell-j} \leq 0 } \right]=: \sum_{j=\nicefrac{d_n}{\alpha}}^\ell L_jF_{\ell-j}.\label{3.2}
\end{align}
Like in the proof of Lemma \ref{lem2.1}, we distinguish cases $0 <\zeta <1 $ and $1  \leq \zeta <2 $.  \\
\textit{When $0<\zeta <1$ } then ${d_n}>A \sqrt j$ for any $A>0$, so applying Lemma \ref{beautiful} like for \eqref{2.3bm} we obtain for ${ d_n \leq \epsilon j}$
\begin{align}
\sum_{j=\nicefrac{d_n}{\alpha}}^\ell L_jF_{\ell-j}& \asymp \sum_{j=\nicefrac{d_n}{\alpha}}^\ell F_{\ell-j} \frac{e^{(\epsilon j \wedge d_n) \Gg( (\epsilon j \wedge d_n)/j)} }{d_n}     \asymp  \frac{1}{\log n}e^{(\log n) \Gg(\log n/ \ell)},  \label{3.1d}
\end{align}
where for any $x$, $g(x):=f(x)-1$. For $ d_n > \epsilon j$, a Markov inequality gives
\begin{align}
\sum_{j=\nicefrac{d_n}{\alpha}}^\ell L_jF_{\ell-j} \leq  C_+ e^{- \epsilon d_n /2} \ell^{-3/2} 
 \leq  C_+ e^{ - \epsilon \log n/ 2 } \ell^{-3/2}, \label{futur1} 
\end{align}
so as $\zeta>0$, considering \eqref{3.1d}
\begin{align}
\sum_{j=\nicefrac{d_n}{\alpha}}^\ell L_jF_{\ell-j} \asymp \frac{e^{(\log n)\Gg(\log n/ \ell)}}{\log n}. \label{zetainf1}
\end{align}
\textit{When $1 \leq \zeta <2$ }, first Lemma \ref{B.1} and \eqref{F} give
\begin{align*}
\sum_{j=\nicefrac{d_n}{\alpha}}^{(1- \epsilon)\ell} L_jF_{\ell-j}\leq {C_+} (\epsilon \ell)^{-3/2} \sum_{j=\nicefrac{d_n}{\alpha}}^{(1- \epsilon)\ell} P(\underline S_j>0) \leq C_+  \ell^{-3/2}\sum_{j=\nicefrac{d_n}{\alpha}}^{(1- \epsilon)\ell} j^{-1/2} \leq C_+  \ell^{-1},
\end{align*}
also
\begin{align}
\sum_{j=(1- \epsilon)\ell}^{\ell} L_jF_{\ell-j}\leq C_+ \sum_{j=(1- \epsilon)\ell}^{\ell} F_{\ell-j}  P(\underline S_j>0) \leq C_+ \sum_{j=(1- \epsilon)\ell}^{\ell}  ( \ell-j+1)^{-3/2} j^{-1/2} \leq C_+  \ell^{-1/2}. \label{3.3a}
\end{align}
For any $j \geq (1-\epsilon) \ell$, $d_n \leq B j^{1/2}$ for all $B>0$, and let $A>B$ 
 \begin{align}
\sum_{j=(1- \epsilon)\ell}^{\ell} L_jF_{\ell-j}&  \geq  \sum_{j=(1- \epsilon)\ell}^{\ell} F_{\ell-j}  P(B j^{1/2} < S_j \leq Aj^{1/2} , \underline S_j>0)  \nonumber\\
 & = \sum_{j=(1- \epsilon)\ell}^{\ell} F_{\ell-j} \sum_{k=B j^{1/2}}^{A j^{1/2}-1}    P(S_j\in (k,k+1] , \underline S_j>0) \label{3.3b}
\end{align}
then Lemma \ref{beautiful} yields
 \begin{align*}
\sum_{k=B j^{1/2}}^{A j^{1/2}-1}    P(S_j\in [k,k+1) , \underline S_j>0) & \asymp \sum_{k=B j^{1/2}}^{A j^{1/2}} k j^{-3/2} e^{- k^2/(2 \sigma^2 j)}   \geq C_- j^{-1/2}.  
\end{align*}
Inserting this in \eqref{3.3b} give $\sum_{j=(1- \epsilon)\ell}^{\ell} L_jF_{\ell-j} \geq C_- \ell^{-1/2}$, this together with \eqref{3.3a} and the above inequality implies
\begin{align}
\sum_{j=\nicefrac{d_n}{\alpha}}^\ell L_jF_{\ell-j}\un_{ d_n \leq \epsilon j} & \asymp \ell^{-1/2} \asymp \ell^{-1/2} e^{(\log n)\Gg(\log n/ \ell)} \label{zetasup1}.
\end{align}
Collecting \eqref{3.2}, \eqref{zetainf1} and \eqref{zetasup1} yields
 \begin{align*}
A_n^+ \leq E \left[e^{S_\ell- \bar S_\ell} \un_{ {\bar S_\ell} > d_n}\right]  \asymp  e^{(\log n)\Gg(\log n/ \ell)} \left( (\log n)^{-1} \un_{0<\zeta <1}+\ell^{-1/2} \un_{1 \leq \zeta <2} \right).
\end{align*} \\
$\bullet$ \textit{ For $A_n^-$}, with $\mathcal{B}_\ell:= \{\sum_{i=1}^\ell e^{S_i}  \leq \ell^\epsilon e^{\bar{S}_\ell} \}$ and $b_n:=\log(c_+ n)$
\begin{align}
\ell^{\epsilon}A_n^- & \geq E \left[{e^{S_\ell-\bar S_\ell}} \un_{\bar{S}_\ell >  b_n,\mathcal{B}_\ell }\right]=E \left[{e^{S_\ell-\bar S_\ell}} \left(\un_{\bar{S}_\ell >  b_n }- \un_{\bar{S}_\ell >  b_n,\overline{\mathcal{B}}_\ell}\right)\right]:=\Gamma_1-\Gamma_2. \label{3.5} 
\end{align} 
$\Gamma_1$ can be treated as $E \left[{e^{S_\ell-\bar S_\ell}} \un_{\bar{S}_\ell >  d_n}\right]$ so  
\begin{align}
\Gamma_1 \asymp e^{(\log n)\Gg(\log n/ \ell)} \left( (\log n)^{-1} \un_{0<\zeta  < 1}+ \ell^{-1/2} \un_{1 \leq \zeta <2} \right) . \label{un}
\end{align}
Recalling that  $\Se_j=\{S_j= \bar S_\ell, \bar S_{j-1}<S_j\}, \, \Gamma_2=\sum_{j=1}^{\ell}   E \left[{e^{S_\ell-S_j}} \un_{S_j >  b_n,\overline {\mathcal{B}}_\ell, \Se_j}\right]$. Note that on $\Se_j$, $\overline{\mathcal B}_\ell=\lbrace Y_1(j)+Y_2(j)> \ell^\epsilon \rbrace$ where $Y_1(j):= \sum_{i=1}^je^{S_i-{S}_j}$ and $Y_2(j):= \sum_{i=j+1}^\ell e^{S_i-{S}_j}$. As for  $0<\delta <1/2$, $\overline{\mathcal B}_\ell\subset\{Y_1(j) > \delta \ell^\epsilon  \} \cup \{Y_2(j) > \delta \ell^\epsilon  \} $
\begin{align}
\Gamma_2&\leq \sum_{j=b_n/ \alpha }^{\ell}E \left[e^{S_\ell-S_j} \un_{{S_j}  >  b_n,\Se_j}\left(  \mathds{1}_{Y_1(j) > \delta \ell^{\epsilon}}+\mathds{1}_{Y_2(j) > \delta \ell^{\epsilon}}\right)\right]=: \sum_{j=b_n/ \alpha }^{\ell} (\Pi_j+\Omega_j).  \label{deux}
\end{align}
\underline{\textit{For} $\sum_{j=b_n/ \alpha }^{\ell} \Omega_j$}, conditioning by $\sigma(S_k,0\leq k\leq j)$  
\begin{align*}
\Omega_j 
=P(\He_j)E\left[e^{S_{\ell-j}} \un_{\bar{S}_{\ell-j} \leq 0,  Y^+(\ell-j) >  \delta \ell^\epsilon  } \right]
\end{align*}
where $Y^{ \pm }(k):= \sum_{i=1}^{k} e^{\pm S_i}$ and $\He_m:=\{ {{S}_m}  >  b_n, \underline S_{m} > 0\}$. Moreover,  using \eqref{F} and the fact that, on $\{\bar{S}_{\ell-j} \leq 0\}$,  $\lbrace Y^+(\ell-j) >  \delta \ell^\epsilon\rbrace=\emptyset$ for $j \leq  \ell-\delta\ell^\epsilon$ 
\begin{align*}
E\left[e^{S_{\ell-j}} \un_{\bar{S}_{\ell-j} \leq 0, Y^+ (\ell-j) >  \delta \ell^\epsilon  } \right ] \leq C_+(\ell-j+1)^{-3/2}\un_{\ell-j \geq \delta \ell^{\epsilon} } , 
\end{align*}
and $P\left( \He_j\right)$ can be treated like $L_j$. Then, for $n$ large enough  
\begin{align}
\sum_{j=b_n/ \alpha }^{\ell} \Omega_j \leq \frac{e^{(\log n)\Gg(\log n/ \ell)}}{\delta \ell^{\epsilon}}  \left( (\log n)^{-1} \un_{0<\zeta <1}+\ell^{-1/2} \un_{1 \leq \zeta <2} \right). \label{3.5b}
\end{align}
\noindent \underline{\textit{For} $\sum_{j=b_n/ \alpha }^{\ell} \Pi_j$}, we have 
$ \Pi_j=P\left({\He_j,  Y^-(j) >  \delta \ell^\epsilon}\right)F_{\ell-j}.$
 Let  $R:= \sup \{0< k \leq j, S_k \leq \gamma \log \ell \}$ with $\gamma>1$, $\tau^+_{ x}:= \inf\{k>0,\ S_k \geq x\}$ for $x\in\mathbb R$   and $t_n=  (\log \ell)^2 $, then
\begin{align*}
\Pi_j =  (P\left({\He_j,  Y^-(j) >  \delta \ell^\epsilon,R\leq \tau^+_{ t_n}}\right)+P\left({\He_j,  Y^-(j) >  \delta \ell^\epsilon,R> \tau^+_{ t_n}}\right))F_{\ell-j}=:\Gamma_3+\Gamma_4.
\end{align*}
\begin{figure}[htbp]
  \begin{center}
    \leavevmode
    \subfloat[First case $(3)$]{%
      \label{Image1}
      \includegraphics[width=6cm]{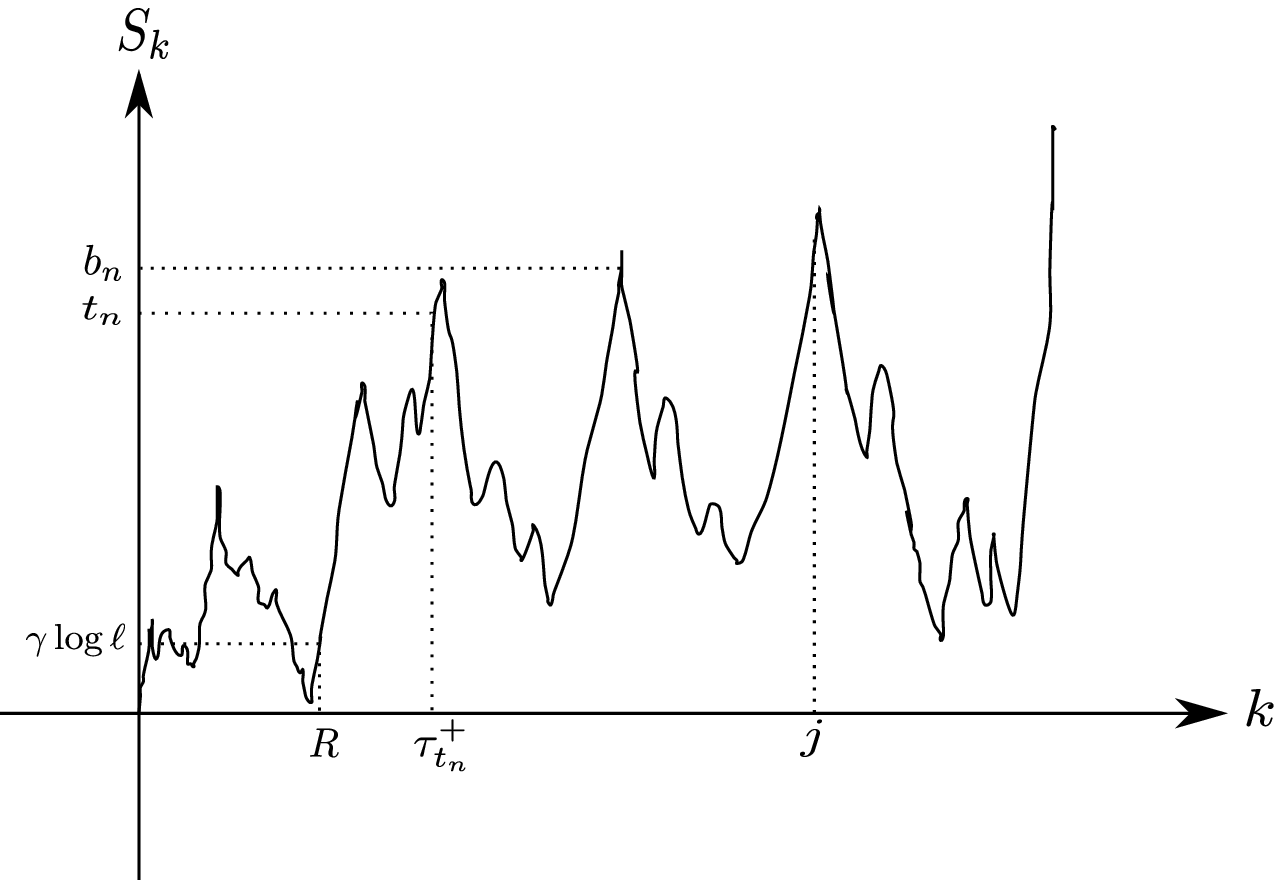}}
    \hspace{2cm}
    \subfloat[Second case $(4)$]{%
      \label{Image2}
      \includegraphics[width=6cm]{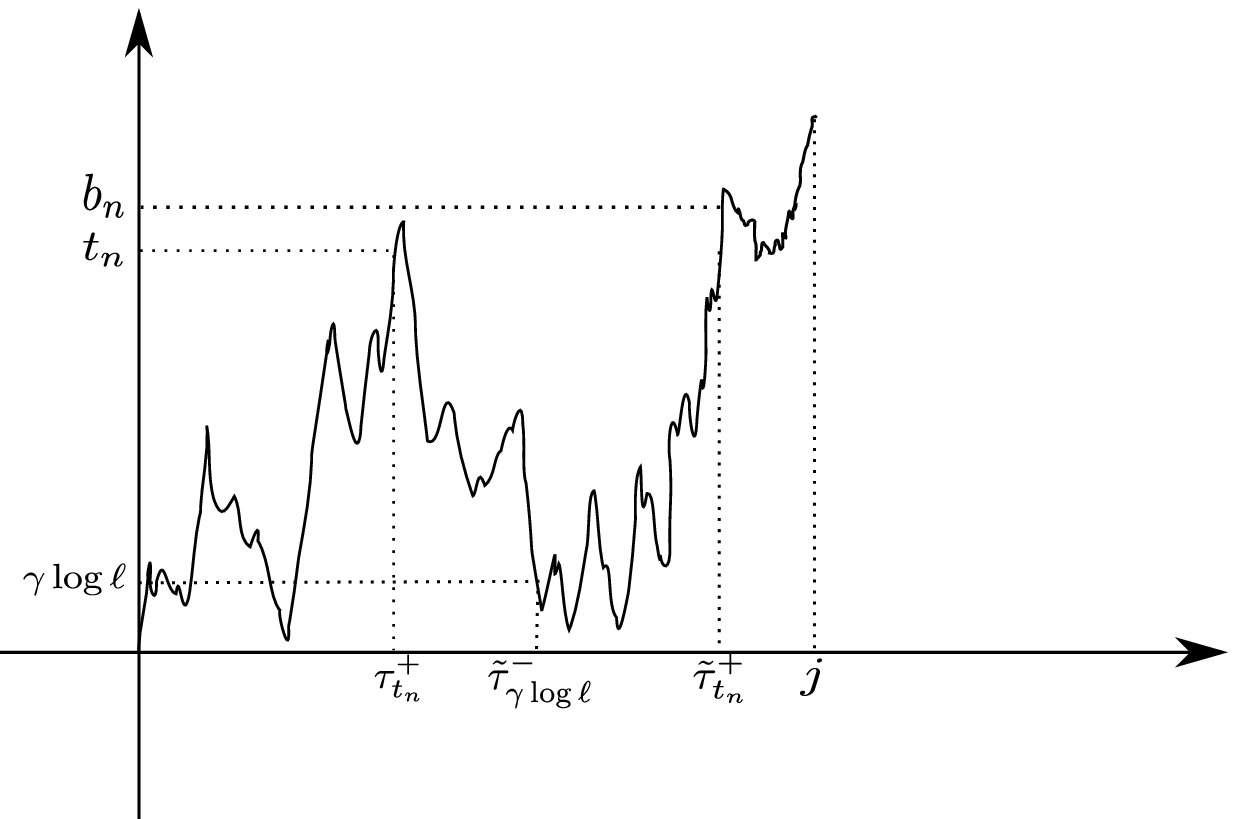}}
          \caption{Two cases}
    \label{Les deux images}
  \end{center}
\end{figure}
\textit{ Upper bound for} $\Gamma_3$, let $\tau^-_{ x}:= \inf\{k>0,\ S_k \leq x\}$ for $x\in\mathbb R$. Notice that on $\{R \leq \tau^+_{ t_n} \}$, $\sum_{k=\tau^+_{t_n}+1}^je^{-S_k}\leq \ell^{1-\gamma }$ implying that   $\{Y^-(j) >  \delta \ell^\epsilon \} \subset \{Y^-(\tau^+_{ t_n})>  \delta \ell^\epsilon/2\} $. 
 Thus, using strong Markov property 
\begin{align}
 &P\left({\He_j, Y^-(j) >  \delta \ell^\epsilon ,R\leq \tau^+_{ t_n}}\right) \leq P\left({\He_j,  Y^-(\tau^+_{ t_n}) >  \delta \ell^\epsilon/2,R \leq \tau^+_{ t_n} }\right) \nonumber \\
&\leq\sum_{k={t_n/\alpha}}^j P\left({\He_j,  Y^-(\tau^+_{ t_n}) >  \delta \ell^\epsilon/2,\tau^+_{ t_n}=k}\right) \nonumber \\
&\leq \sum_{k={t_n/\alpha}}^j\sup_{z \in [0,\alpha]}P_{z+t_n}\left(\He_{j-k}\right) P({ Y^-(\tau^+_{ t_n}) >  \delta \ell^\epsilon/2,\tau^-_0>\tau^+_{ t_n}=k}). \label{3.6bm}
\end{align}
Case $\zeta <1$, we use the following upper bound for  \eqref{3.6bm}  
$$ P(Y^-(\tau^+_{ t_n}) >  {\delta \ell^\epsilon}/{2}, \tau^-_0 > \tau^+_{ t_n} ) \sup_{t_n/\alpha\leq k \leq j} \sup_{z\in[0,\alpha]}P_{z+t_n}\left(\He_{j-k}\right).$$
Lemma \ref{lemshady} gives an upper bound for the first probability. Moreover with the help of Lemma \ref{beautiful2} and a similar reasoning as for \eqref{3.1d} and \eqref{futur1}, for all $t_n/\alpha\leq k \leq j $  
\begin{align*}
 \sup_{z\in[0,\alpha]}P_{z+t_n}\left(\He_{j-k}\right) \leq C_+ { t_n}{b_n}^{-1}e^{(\epsilon j \wedge b_n) \Gg((\epsilon j \wedge b_n)/j)}. 
\end{align*}
So $
 \Gamma_3  \leq C_+t_n b_n^{-1}\ell^{-\epsilon} (\ell-j+1)^{-3/2}e^{(\epsilon j \wedge b_n)\Gg((\epsilon j \wedge b_n)/j)} $.\\
Case $\zeta \geq 1$, here the following upper bound for \eqref{3.6bm} is useful 
\begin{align*}
& P(Y^-(\tau^+_{ t_n}) >  {\delta \ell^\epsilon}/{2}, \tau^-_0 > \tau^+_{ t_n} ) \sup_{{t_n}/{\alpha}\leq k \leq j/2} \sup_{z\in[0,\alpha]}P_{z+t_n}\left(\He_{j-k}\right)+P(\tau^-_0 \wedge \tau^+_{ t_n} \geq j/2) \\
& \leq P(Y^-(\tau^+_{ t_n}) >  {\delta \ell^\epsilon}/{2}, \tau^-_0 > \tau^+_{ t_n} )P_{ \alpha+t_n}( \underline{S}_{j/2} >0)+2E[\tau^-_0 \wedge \tau^+_{ t_n}]/j  \leq C_+t_n\left(\frac{1}{\ell^{\epsilon} j^{1/2}}+\frac{1}{ j}\right),
\end{align*}
with Lemmata \ref{B.1} and \ref{lemshady}.  So
\begin{align*}
\Gamma_3 \leq   \frac{C_+t_n}{(\ell-j+1)^{3/2}} \left(\frac{1}{\ell^{\epsilon} j^{1/2}}+\frac{1}{ j} \right) .
\end{align*}
\textit{Upper bound for} $\Gamma_4$,  first note that on $\{R> \tau^+_{ t_n}\}$, the following hitting times $\tilde \tau^-_{ \gamma \log \ell }:=\inf\{ k\in\rrbracket \tau^+_{t_n},j\rrbracket, S_k \leq \gamma \log \ell \}=\tau^+_{t_n}+\theta_{\tau^+_{t_n}}\circ \tau^-_{ \gamma \log \ell }$, and  $\tilde \tau^+_{t_n }:=\inf\{ k\in\rrbracket \tilde\tau_{\gamma \log \ell }^-,j\rrbracket, S_k \geq t_n   \}=\tilde \tau_{\gamma \log \ell }^-+\theta_{\tau_{\gamma \log \ell }^-}\circ\tau^+_{t_n}$ (where $\theta$ is the shift operator) exist. With these notations according to Lemma \ref{B.1}
\begin{align*}
&P(\tau^-_0>\tilde \tau^+_{t_n})=P(\tau^+_{t_n}<\tau^+_{t_n}+\theta_{\tau^+_{t_n}}\circ \tau^-_{ \gamma \log \ell }<\tilde \tau_{\gamma \log \ell }^-+\theta_{\tilde \tau_{\gamma \log \ell }^-}\circ\tau^+_{t_n}<\tau^-_0)\\
&\leq P(\tau^-_0>\tau^+_{t_n})\sup_{z\in[0,\alpha]}P_{\gamma\log \ell-z}(\tau^-_0>\tau^+_{t_n})\leq C_+\frac{\gamma\log\ell}{(t_n)^2}.
\end{align*}
Again at this point we distinguish the cases $\zeta <1$ or $\zeta \geq 1$.

\noindent When $\zeta <1$ the above inequality yields 
\begin{align*}
P\left({\He_j,R> \tau^+_{ t_n}}\right)
&\leq \sup_{{2t_n}/{\alpha}\leq k \leq j} \sup_{z\in[0,\alpha]}P_{z+t_n}\left(\He_{j-k}\right)P(\tau^-_0>\tilde \tau^+_{t_n})\\
& \leq C_+\frac{\gamma\log\ell}{(t_n)^2} \sup_{{2t_n}/{\alpha}\leq k \leq j} \sup_{z\in[0,\alpha]}P_{z+t_n}\left(\He_{j-k}\right).
\end{align*}
Finally using that $P\left({\He_j,  Y^-(j) >  \delta \ell^\epsilon,R> \tau^+_{ t_n}}\right) \leq P\left(\He_j,R> \tau^+_{ t_n}\right)$, Lemma \ref{beautiful2} and \eqref{F} 
\begin{align*}
\Gamma_4 &
\leq C_+ \frac{e^{(\epsilon j \wedge b_n)\Gg((\epsilon j \wedge b_n)/j)}}{b_n ( \ell-j+1)^{3/2}\log \ell}. 
\end{align*}
 When $\zeta \geq 1$, we have 
 $$P\left({\He_j,  Y^-(j) >  \delta \ell^\epsilon,R> \tau^+_{ t_n}}\right)  \leq P\left({\He_j,j/2 \geq R> \tau^+_{ t_n}}\right)+P\left({\He_j , R> \tau^+_{ t_n}} \vee j/2 \right). $$ 
 Moreover $P\left({\He_j,j/2 \geq R> \tau^+_{ t_n}}\right) \leq P\left( \tau^+_{ t_n}< \tau^-_0\right)\sup_{-\alpha\leq  x \leq \alpha}P_{\gamma \log \ell+x } (\underline{S}_{j/2}>0)$ \\
 and $P\left({\He_j, R> \tau^+_{ t_n}} \vee j/2 \right) \leq P(\underline{S}_{j/2}>0)\sup_{-\alpha\leq  x \leq \alpha}P_{\gamma \log \ell+x }( \tau^-_0>\tau^+_{b_n})$. So using Lemma \ref{B.1} 
  \begin{align*}
\Gamma_4 &
\leq \frac{ C_+}{j^{1/2}(\log \ell )( \ell-j+1)^{3/2}}. 
\end{align*}
Collecting \eqref{3.5b} and $\Gamma_3$, $\Gamma_4$ with \eqref{deux}
\begin{align*}
\Gamma_2 \leq  C_+   \frac{e^{b_n \Gg(b_n/\ell)}}{ b_n \log \ell} \sim C_+ \frac{e^{(\log n)\Gg(\log n/\ell) }}{ (\log n) \log \ell },
\end{align*}
if $\zeta <1$ and 
 \begin{align*}
\Gamma_2 \leq     \frac{C_+} {\ell^{1/2}\log \ell} \leq  C_+ \frac{e^{(\log n)\Gg(\log n/\ell)  }}{ \ell^{1/2}\log \ell },
\end{align*}
if $ \zeta \geq 1$. $\Gamma_2$ is therefore negligible compared to $\Gamma_1$ (see \eqref{un}) and \eqref{3.5} implies that $A^-_n \geq C_+ e^{(\log n)\Gg((\log n)^{- \zeta}) }\left((\log n)^{-1-\varepsilon} \un_{0<\zeta <1}+ \ell^{-1/2-\varepsilon}\un_{1 \leq \zeta <2} \right) $. This with Lemma \ref{ultimlem} finish the proof of \eqref{1.7}.

\subsection{From   $\boldsymbol{ \Ks_{\Phi(n)}(\ell)}$ to $\boldsymbol{K_n(\ell)}$   and $\boldsymbol{\Kt_n(\ell)}$ (proof of (\ref{1.6}) and (\ref{1.8}))}


Let $\Phi_1(n):=(1-2\epsilon) \log n$, we need the following

\begin{Lem} \label{lem3.1}  Let $\mathcal{A}:= \{|z|=\ell, \overline{V}(z) \leq \Phi_1(n) \}$,
$$\lim_{n \rightarrow + \infty} \p\left(\min _{z \in \A} \lo(z, T_{\phi}^{n^{1- \epsilon}}) \geq 1  \right)=1,$$
 which implies $\lim_{n \rightarrow + \infty } \p(  K_{n^{1- \epsilon}}(\ell) \geq \Ks_{\Phi_1(n)}(\ell))=1.$
\end{Lem}

\begin{Pre}
Applying Corollary \ref{routine}, 
\begin{eqnarray*}
\pe(\cup _{z \in \A}\{ \lo(z, T_{\phi}^{n^{1- \epsilon}})=0 \} )\leq \vert  \mathcal{A}\vert e^{-c_- n^{\epsilon} \ell^{-1}} \leq \Ks_{\Phi_1(n)}(\ell) e^{-c_- n^{\epsilon/2}}
 \end{eqnarray*}
Using \eqref{Prop2}, $E[\Ks_{\Phi_1(n)}(\ell)] \leq e^{\Phi_1(n)}$ and the proof is achieved.
\end{Pre} \\
The above Lemma together with \eqref{2.7lim} (taking $\Phi(n)=\Phi_1(n)$),  give for $n$ large enough
\begin{align}
\lim_{n \rightarrow + \infty} \p\left( \max_{|z|=\ell- \log n/\tilde{\gamma} } \min_{ y \in \mathcal{C}_{\ell} (z)} {\lo(y,T_{\phi}^{n^{1- \epsilon}}) \geq 1}   \right) =1
\end{align} 

\noindent To obtain the lower bound in \eqref{1.8} we finally use the following result that can be deduced from \cite{HuShi10b}  (see \cite{AndreolettiDebs1} Lemma 3.2 and what follows for details)
\begin{align}
\forall \delta>0,\,  \lim_{n \rightarrow + \infty}\p( \lo(\phi,n) \geq n^{1-\delta} )=1.
\label{lem3.2}
\end{align} 
For the lower bound in \eqref{1.6}, we use Lemma \ref{lem3.1}, \eqref{lem3.2} and finally the lower bound in \eqref{Prop2}.

For the upper bound in \eqref{1.6}, denote $u_n:=C \left( {\log \log n} \vee (\log n)^{1-\zeta}\right)$, where $C>0$. As $n \leq T_{\phi}^n$, by Markov inequality and \eqref{1.7},  $\p \left(  \log \Kt_{n}(\ell)   \geq  \log n-u_n \right) \leq \p( K_n(\ell)  \geq n e^{-u_n }  ) \leq C_+ e^{u_n } e^{ -(\log n)^{1-\zeta}/2 \sigma^2}\ell^{-1/2} $ which gives the upper bound adjusting $C$ properly.

\section{Visited points along the GW}

In this paragraph we study the manner the random walk visits  the tree.
  
\subsection{Visits of clusters at deterministic cuts (proof of (\ref{1.10b}))}
\begin{figure}[h!]
\begin{center}
\includegraphics[width=7cm]{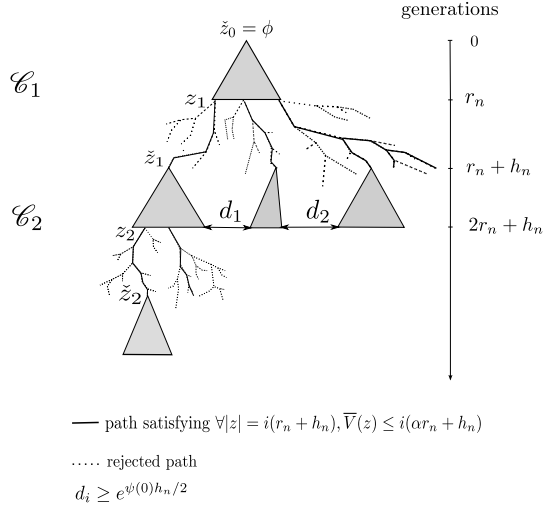}
\end{center}
\caption{Clusters at regular cuts} \label{fig3}
\end{figure}

\begin{figure}[h!]
\begin{center}
\includegraphics[width=11cm]{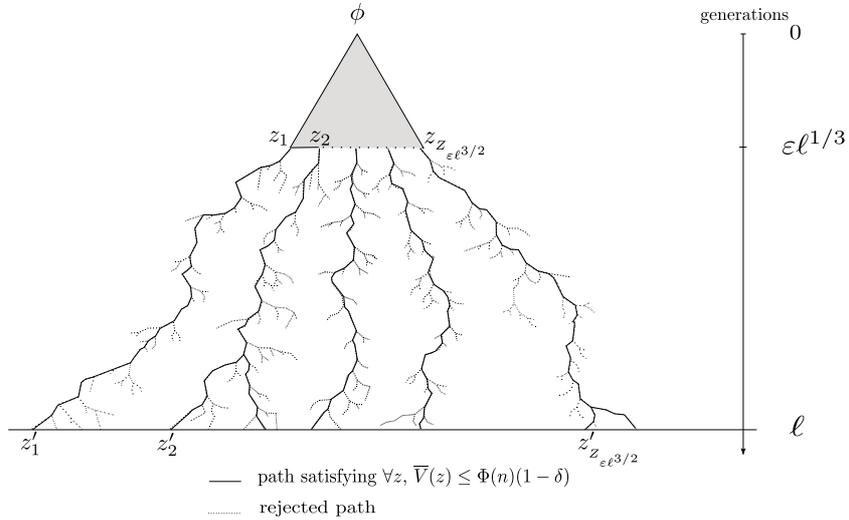}
\end{center}
\caption{Distant visited sites} \label{fig4}
\end{figure}

Recall that a cluster initiated at $z$ with end generation $m$ is the set  $\Sl_m(z)=\lbrace u>z,\vert u\vert=m \rbrace$. Also take $k_n = \Phi(n)^{{\bf k }}$, $r_n = \Phi(n)^{{\bf r }}$, $s_n = \Phi(n)^{{\bf s }}$ with  ${\bf s}>0$ and  $h_n$  sequences such that
\begin{align} 
h_n=\frac{\ell-k_nr_n}{k_n-1},\  k_n(\alpha r_n+s_n)-s_n \leq \Phi(n)(1- 2\epsilon), \label{4.1b}
\end{align} 
where $\alpha:=\vert \log \varepsilon_0\vert$ (see \eqref{remA1} for details).\\
 We define recursively clusters at generations $ir_n+(i-1)h_n $ for all $1\leq i \leq k_n$ in the following way (see also Figure \ref{fig3}): the iteration starts with $\check{z}_0=\phi$ and
\begin{align*}
\forall z_i\in\Sl_{ir_n+(i-1)h_n}(\check{z}_{i-1}), \check{z}_i={\bf inf}\lbrace u>z_i, \vert u\vert=i(r_n+h_n),  \,\overline{V}(u) \leq i(\alpha r_n+s_n)\rbrace.
\end{align*}

The individuals of these clusters form a subtree of the GW, moreover for all $z$  of this subtree at generation $\ell$, $\overline{V}(z) \leq \Phi(n)(1- 2\epsilon)$.
For a fixed $i\in\llbracket 1 ; k_n\rrbracket$,  $\Cs_{i}$ denotes, among the previously defined clusters,  the ones rooted at generation $(i-1)(r_n+h_n)$ and with end points at generation $i r_n+(i-1)h_n$, in other words $\Sl_{i r_n+(i-1)h_n}(.)$. We first give an upper bound for the probability that for all $i \leq k_n$ every clusters in $\Cs_{i}$ are fully visited before $T_{\phi}^{n^{1- \epsilon}}$
\begin{align}
 \pe\left(\bigcup_{i=1}^{k_n}\bigcup_{\Dl \in \Cs_{i} }  \bigcup_{z\in\Dl}\left\{  \lo(z,T_{\phi}^{n^{1- \epsilon}})=0\right\}  \right) 
& \leq \sum_{i=1}^{k_n}\sum_{\Dl \in \Cs_{i} }   \sum_{z \in \Dl }   \pe\left( \lo(z,T_{\phi}^{n^{1- \epsilon}}) =0 \right).  \nonumber 
\end{align}
In the previous formula $\bigcup_{i=1}^{k_n}\bigcup_{\Dl \in \Cs_{i} }$ is an abuse of notation as the sets of clusters are defined recursively. 
With a similar reasoning as the one for Corollary \ref{routine} and the ellipticity condition for the number of descendants 
\begin{align}
  \pe\left(\bigcup_{i=1}^{k_n}\bigcup_{\Dl \in \Cs_{i} }  \bigcup_{z\in\Dl}\left\{  \lo(z,T_{\phi}^{n^{1- \epsilon}})=0\right\}  \right) & \leq \exp(k_nr_n\log N_0 -c_- n^{1- \epsilon} e^{-\Phi(n)(1-2\epsilon)}/\ell ).  \label{4.0}
\end{align}

\noindent We now prove the existence of such clusters, this implies new  constraints on ${\bf k}$, ${\bf r}$ and ${\bf s}$ in addition to \eqref{4.1b}. 


\noindent First, ellipticity conditions imply that for any site $z>y$,  $ {V}(z)-V(y) \leq \alpha(\vert z\vert-\vert y\vert) $ a.s. 
Thus, for all $a\in \mathbb N$ and $b>0$,  $\lbrace \exists \check{z}>z,\vert \check{z}\vert=\vert z\vert +a,\tilde {V}(z,\check{z})\leq b\rbrace $ is a.s. contained in $\lbrace {V}(z)-V(y)\leq \alpha(\vert z\vert-\vert y\vert), \exists \check{z}>z,\vert \check{z}\vert=\vert z\vert +a,\tilde{V}(y,\check{z})\leq \alpha(\vert z\vert-\vert y\vert)+b\rbrace$ ($\tilde {V}$ is defined in the proof of \eqref{2.9}).\\
Then, with our slight abuse of notation, a.s. $\Aa:=\{$the clusters $\Sl_.(.)$ exist$\}$ contains  
\begin{align*}
\bigcap_{i=1}^{k_n}\left\lbrace\bigcap_{z_i\in\mathcal C_{ir_n+(i-1)h_n}(\check{z}_{i-1})}\left\lbrace \exists \check{z}_i>z_i,\vert \check{z}_i\vert=i(r_n+h_n), \tilde{V}(z_i,\check{z}_i)\leq s_n \right\rbrace\right\rbrace.
\end{align*}
The independence of the increments of $V$ and the ellipticity assumptions on the number of descendants ($|\Sl_.(.)| \leq N_0 ^{r_n}  $) imply
\begin{align*}
P(\Aa)& \geq \left[P(\exists |{z}|=h_n, \overline{V}({z}) \leq  s_n )\right]^{N_0 ^{k_n r_n}}
=\left [P\left(\min_{|{z}|=h_n}\overline{V}({z}) \leq  s_n \right)\right]^{N_0 ^{k_n r_n}}.  
\end{align*}
Assuming $ {\bf s } \geq (1+\zeta-{\bf k })/3$, Lemma \ref{ProD3} yields $P(\Aa)\geq  (1-\lambda_{h_n})^{N_0 ^{k_n r_n}}\sim \exp\left( -N_0 ^{k_n r_n}\lambda_{h_n}  \right)$. To choose $\bf k$, $\bf r$ and $\bf s$, we have to take into account the last constraint in \eqref{4.1b}, $\bf s\leq 1- \bf k$ and that if $q_0+q_1\neq 0$, $\bf k+\bf r<s$. We distinguish two cases 
 \begin{itemize}
 \item if $ 0 < \zeta  \leq 1 $, let $0<\delta<\zeta/2 $, take $ {\bf s} =(1+ \zeta)/2- \delta$, ${\bf k} =\delta /2$, and ${\bf r} ={(1-\zeta)}/{2}+ \nicefrac{\delta}{2}$,
\item if $1< \zeta  < 2$, let $0<\delta<(2-\zeta)/3$, take ${\bf s}=(1+\zeta)/3$, ${\bf k}=\delta$ and ${\bf r}=(1+\zeta-4 \delta)/3$.
\end{itemize}
Thus in both cases
\begin{align}\label{bc}
P(\Aa) &\geq 1 -C_+ e^{-c_1 s_n}N_0 ^{k_n r_n} \underset{n\rightarrow+\infty}{\longrightarrow} 1.  
\end{align} 
When $q_0+q_1 = 0$, the above choices give an even better rate of convergence for $P(\Aa)$.  \\
We now move back to \eqref{4.0}, \eqref{bc} together with the fact that   $\Phi(n) \leq \log n+o(\log n)$  implies   
\begin{align*}
\lim_{n \rightarrow + \infty} \p \left(\bigcap_{i=1}^{k_n}\bigcap_{\Dl \in \Cs_{i} }\bigcap_{z \in \Dl} \left\{ {\lo(z,T_{\phi}^{n^{1- \epsilon}}) \geq 1}  \right\} , \Aa \right) = 1.  
\end{align*} 
According to \eqref{lem3.2},  as $\p(n \geq T_{\phi}^{n^{1- \epsilon}})$ tends to one we finally obtain
 \begin{align*}
\lim_{n \rightarrow + \infty} \p \left(\bigcap_{i=1}^{k_n}\bigcap_{\Dl \in \Cs_{i} } \bigcap_{z \in \Dl }\left\{ {\lo(z,n) \geq 1} \right\} , \Aa \right) = 1.
\end{align*}
So we can find set of clusters at regular cuts on the tree which are fully visited. To finish the proof of \eqref{1.10b} we first show the existence of a lower bound for the number of visited clusters. Using successively that conditionally on $\Aa,\, |\Cs_{i}|$ is equal in law to  $Z_{(i-1)r_n}$, Theorem \ref{ND} and \eqref{bc} 
\begin{align*} 
 P\left(\exists i\in\llbracket 2, k_n\rrbracket, |\Cs_{i}|\leq e^{\psi(0) (i-1) r_n/2 }  \right)& \leq  P\left(\exists i\in\llbracket 2, k_n\rrbracket, |\Cs_{i}|\leq e^{\psi(0) (i-1) r_n/2 } ,   \Aa \right)+P( \overline{\Aa}) \\
  & \leq  \sum_{i=2}^{k_n} P\left(Z_{(i-1) r_n} \leq e^{\psi(0) (i-1) r_n/2 } \right)+P( \overline{\Aa})\\
  &\leq e^{-\psi(0) \nu r_n/4 }+C_+ e^{-c_1 s_n} N_0 ^{k_n r_n}. 
\end{align*}
Note that for the first term we have used the left tail of $Z_.$ with $q_0+q_1>0$ as the other case provide an even better rate of convergence.   
Finally we prove that the previously defined visited clusters are very spaced out. Recalling the definition of ${\bf D}$ before Theorem \ref{Th2}, 
\begin{align*}
   P\left(\exists i\in\llbracket 2, k_n\rrbracket, {\bf{D}}(\Cs_{i})  
  \leq e^{\psi(0) h_n/2 } \right) 
  &\leq  P( \overline{\Aa})+\sum_{i=2}^{k_n} P\left({\bf{D}}(\Cs_{i}) \leq e^{\psi(0) h_n/2 },\Aa  \right) \\
 \leq P( \overline{\Aa})&+\sum_{i=2}^{k_n} P\left(\bigcup_{\Dl \in \Cs_{i-1} } \bigcup_{z \in \D } \{ |\mathcal{C}_{i (r_n+h_n)}(z)  | \leq e^{\psi(0) h_n/2 } \}\cap \Aa  \right).
\end{align*}
As conditionally on $\Aa,\,  |\mathcal{C}_{i( r_n+ h_n)}(z)  |$ and $Z_{h_n}$ are equal in law, on  $\mathscr{D}_i:=\{ |\Cs_{i}|\leq e^{2 \psi(0) (i-1) r_n }\},$ $|\mathcal D| \, |\Cs_i|\leq N_0^{r_n}e^{2 \psi(0) (i-1) r_n }$ so Theorem \ref{ND} yields 
\begin{align*} 
P\left({\bf{D}}(\Cs_{i}) \leq e^{\psi(0) h_n/2 },\Aa  \right) 
& \leq P\left(\bigcup_{\Dl \in \Cs_{i-1} } \bigcup_{z \in \D } \{ |\mathcal{C}_{i( r_n+ h_n)}(z)  | \leq e^{\psi(0) h_n/2 } \}\cap \Aa \cap\mathscr{D}_i  \right)+P(\overline{\mathscr{D}}_i) \\
& \leq e^{r_n \log N_0 +2 \psi(0) (i-1) r_n } P( Z_{h_n} \leq e^{\psi(0) h_n/2 } )+ e^{- \psi(0) (i-1) r_n } \\
& \leq e^{r_n \log N_0 +2 \psi(0) (i-1) r_n - \nu \psi(0) h_n/2 } + e^{- \psi(0) (i-1) r_n }.
\end{align*}
Consequently
$$P\left(\exists i\in\llbracket 2, k_n\rrbracket, {\bf{D}}(\Cs_{i})  
  \leq e^{\psi(0) h_n/2 }   \right) \leq e^{3 \psi(0) k_n r_n - \nu \psi(0) h_n/2 }+2e^{-\psi(0)r_n}+C_+e^{-c_1 s_n}N_0^{k_n r_n},$$
  moreover as $h_n \sim (\log n)^{1+ \zeta-\bf{k}} $, ${\bf k}+{\bf r} < {\bf s}<1$ and ${\bf k} < \zeta$ we obtain the result.

\subsection{Proof of (\ref{1.9}) }
Let $m=\epsilon \ell^{1/3} $, $\delta>0$, define $\B$ the set of points $z'$ such that for all $|z|=m $, $z':= {\bf inf}\{u>z, |u|=\ell,  \overline{V}(u) \leq \Phi(n)(1- \delta)\}$. Corollary \ref{routine} gives
$$
\pe\left(\bigcup_{z' \in \B} \left\{ \lo(z^\prime,T_{\phi}^{n^{1- \epsilon}}) = 0 \right\} \right)  
\leq |\B|e^{-c_- n^{1- \epsilon} e^{-\Phi(n)(1-{\delta})}/\ell}\leq Z_me^{-c_- n^{1- \epsilon} e^{-\Phi(n)(1-{\delta})}/\ell} . 
$$
As $\Phi(n) \leq \log n + o(\log n)$ and $E[Z_m]=e^{\psi(0)m}$, taking $\delta=2 \epsilon$ 
\begin{align}
\p\left(\bigcup_{z' \in \B} \left\{ \lo(z^\prime,T_{\phi}^{n^{1- \epsilon}})= 0 \right\} \right)  \leq e^{-c_- n^{\epsilon/2}}\label{4.7}.
\end{align}
We now prove that $\lim_{n \rightarrow + \infty} P(|\B| =Z_m)=1$.
From \cite{McDiarmid} (see also \cite{AndreolettiDebs1} Lemma 2.1), $\lim_{n \rightarrow + \infty}P( \max_{|z|=m}\overline{V}(z) \leq 2 \tilde \gamma m )=1$, so as for $n$ large enough $\ell^{1/3}/\Phi(n) \leq \delta$ with the same arguments used in the proof of Lemma  \ref{ProD3} 
\begin{align*}
P(|\B|<Z_{m})  &\leq  P\left( \bigcup_{ |z|=m} \left\{ \forall z'>z, |z'|=\ell,\overline{V}(z') > (1- \delta)\Phi(n) \right\} \right) \\
& \leq 1-\left(1-P\left( \min_{|z|=\ell-m}\overline{V}(z) > (1- 2\delta)\Phi(n) \right)\right)^{e^{2  \psi(0) m}}  +e^{-\psi(0) m}\\
& \leq 2 P\left( \min_{|z|=\ell-m} \overline{V}(z) > (1-  4\epsilon) \Phi(n) \right) e^{2  \psi(0) m } +e^{-\psi(0) m}.
\end{align*}
To finish we put ourself in the case $q_0+q_1>0$ (the other case is treated similarly), using Lemma \ref{ProD3}
$$ P\left( \min_{|z|=\ell-m}\overline{V}(z) > (1- 4\epsilon) \Phi(n)\right) \leq e^{- c_1 \ell^{1/3} }. $$ 
 We can now choose $\epsilon$ small enough and obtain, $P(|\B|=Z_{m}) \geq 1-e^{-c_1 \ell^{1/3}/2 }  $. Moving back to \eqref{4.7} $\p( \forall z \in \B, \lo(z,T_{\phi}^{n^{1- \epsilon}})  \geq 1 , |\B|=Z_m ) \geq 1-o(1)$.
Finally to obtain \eqref{1.9} we apply \eqref{lem3.2}. 

\appendix

\section{Basic facts for branching processes and Galton-Watson trees}

\subsection{Biggins-Kypriaou identities and properties of $\boldsymbol{\psi}$ \label{apA}}


For any $n\geq1$ and any mesurable function $F:\mathbb R^n\rightarrow[0,+\infty)$, Biggins-Kyprianou identity is given by

\begin{equation}\label{B1}
E\left[\sum_{\vert x\vert =n}e^{-V(x)}F(V(x_i),1\leq i\leq n)\right]=E[F(S_i,1\leq i\leq n)]
\end{equation}
where $(S_i-S_{i-1})_{i\geq 1}$ are i.i.d. random variables, and the law of $S_1$ is determined by
\begin{equation}\label{B2}
E[f(S_1)]=E\left[\sum_{\vert x\vert =1}e^{-V(x)}f(V(x))\right], 
\end{equation}
for any measurable function $f:\mathbb R\rightarrow[0,+\infty)$. A proof can be found in \cite{BigginsKyprianou}, see also \cite{Shi4}. We have the following identities
\begin{align*}
 \psi(t)= \log E[e^{(1-t)S_1}],  \psi'(t)= -\frac{E[S_1e^{(1-t)S_1}]}{E[e^{(1-t)S_1}]}.
\end{align*}
In particular, $E[N]=e^{\psi(0)}=E[e^{S_1}]$ 
and the hypothesis $\psi'(1)=0$ equates to $E[S_1]=0$.
\noindent 
\begin{Rem}\label{remA1}
Let $\alpha:=\vert \log \varepsilon_0\vert$, we have $\p(\vert S_1\vert  \leq \alpha)=1$. Indeed taking $f(x)=\mathds{1}_{\vert x\vert  \leq \alpha}$  and using Biggins-Kiprianou, $\p(\vert S_1\vert \leq\alpha)=E[f(S_1)]=e^{\psi(1)}=1$.
\end{Rem}



\subsection{Left tail of $\boldsymbol{Z_n}$}

Recall that the positive martingale $(W_n)_{n\geq0}:=(\nicefrac{Z_n}{e^{\psi(0)n}})_{n\geq0}$ converges a.s. to a non degenerate limit $W$ (see for instance \cite{Shi4}). Moreover $W$ has a positive continuous density function denoted $w$.
Bingham  \cite{Bingham}  shows that for the Schr\"oder case ($q_0+q_1>0$), there exists $0<\nu<1$ such that for small $x$, $ w(x) \sim x^{\nu-1}$ and
for the B\"ottcher  case ($q_0+q_1=0$) there exists $\beta\in(0,1)$ such that when $x\rightarrow0$,  $\log w(x) \sim -{x^{-\frac{\beta}{1-\beta}}}$. 
The results of \cite{AtNe} and then \cite{Fleiwach2} (Theorems 4 and 5) and \cite{Fleiwach} (Theorem 7) lead to
\begin{The} \label{ND}
Let $0<\kappa < 1$ then
$P(Z_n \leq e^{ \kappa \psi(0) n}) \sim e^{-\nu \psi(0) (1-\kappa)n} $ in the Schr\"oder case, and 
  $\log P(Z_n \leq e^{ \kappa \psi(0) n})\sim\log w\left({e^{\psi(0)(\kappa-1)n }}\right)$ in the B\"ottcher case.
\end{The}


\section{Results for sums of i.i.d. random variables }

In this section we recall basic facts for sum of i.i.d. random variables applied to $(S_n)_{n\geq0}$ of Section A. Recall that for all $x\in\mathbb R$, $\tau_x^+=\inf\lbrace n\geq1,S_ n\geq x\rbrace$ and $\tau_x^-=\inf\lbrace n\geq1,S_ n\leq x\rbrace$. 
The following results are standard and can be found in \cite{Aidekon} and \cite{Spi}. 
  
 \begin{lem} \label{B.1}
 For all $x \in [0,y]$ and $m$ large enough 
\begin{align*}
 P_x(\tau^+_y<\tau^-_0) \asymp \frac{x+1}{y+1},\   E[\tau^+_y \wedge \tau^-_0] \asymp  y\mbox{ and }P_x(\tau^-_0>m) \asymp \frac{ x+1 }{\sqrt m}.
 \end{align*}
 \end{lem}
Recalling that for all $n\geq1$,  $Y^{ - }(n)= \sum_{i=1}^{n} e^{-S_i}$, we have 
\begin{lem} \label{lemshady} There exists  a constant $C_+>1$ such that for all $a\geq 0$ and $M>0$ 
\begin{align*}
 P\left(Y^-(\tau_a^+) > M,\ \tau^+_a<\tau^-_0 \right) \leq \nicefrac{C_+}{M}. 
 \end{align*}
\end{lem}

\begin{Pre}
According to \cite{AiHuZi} p.19, there exists $C_+>1$ such that for all $0\leq a\leq L\leq 1$, 
$$\E\left[Y^-({\tau^-_0\wedge \tau^+_a})\right]\leq C_+\nicefrac{a+1}{a}.$$
As $P\left( Y^-(\tau_a^+) > M,\ \tau^+_a<\tau^-_0 \right)\leq P\left( Y^-({\tau^-_0\wedge \tau^+_a}) > M\right) $, we conclude using the Markov inequality.  
\end{Pre}

\begin{lem} \label{B.3}
 For any $m \geq 1$,
\begin{align}
    E[e^{S_m} \un_{\tau_0^->m}] \asymp  {(m+1)^{-3/2}}.\label{F}
 \end{align}
 \end{lem}
 
 \begin{Pre}
The upper bound can be found in \cite{HuShi24} p.44, the lower bound can be obtained as follows: $E[e^{S_m} \un_{\tau_0^->m}]\leq \sum_{k=0}^{\alpha m} e^{-k}P(S_m \in [-k,-k+1], \bar S_m<0) =\sum_{k=0}^{\alpha m} ke^{-k} m^{-3/2} \asymp   {(m+1)^{-3/2}}$.
\end{Pre}
 
Recalling that for all $n\geq1$,  $Y^{ - }(n)= \sum_{i=1}^{n} e^{-S_i}$, we have 
\begin{lem} \label{lemshady} There exists  a constant $C_+>1$ such that for all $a\geq 0$ and $M>0$ 
\begin{align*}
 P\left(Y^-(\tau_a^+) > M,\ \tau^+_a<\tau^-_0 \right) \leq \nicefrac{C_+}{M}. 
 \end{align*}
\end{lem}

\noindent  The following Lemma may be found in the literature, however as we can prove it easily for our case we present a short proof.
\begin{lem} \label{beautiful2} Let $m>1$, assume that $b=b(m) \geq \sigma^2 \sqrt m \log m$, with $\lim_{m \rightarrow +\infty} b/m=0 $, and $a=a(m)>0$ is such that $ \lim_{m \rightarrow + \infty} a/ \sqrt{m}=0$, then for all $m$ large enough
\begin{align}
P_a\left({S}_m >b, \underline{S}_m >0 \right)    \leq C_+ \frac{a}{b} e^{b \cdot  \Gg(b/m)}, \label{1.5}
\end{align}
with $g(x)=f(x)-1=-\frac{x}{2 \sigma^2}+x^2 \lambda(x).$ For all $\epsilon>0$ and $r> \epsilon m$.
\begin{align}
P_a\left({S}_m \geq r, \underline{S}_m >0 \right) \leq C_+ \frac{a  }{ \sqrt{m}} e^{-sr+m\psi(1-s) }. \label{beautiful3}
\end{align}
 \end{lem}

\begin{Pre}
 For \eqref{1.5}, 
  let $\omega$ a positive function of $b$ and $m$, such that $\omega \leq 2\sqrt m /b$ and that we choose later,  write $P_a\left({S}_m >b, \underline{S}_m >0 \right)$ as
\begin{align*}
P_a\left({S}_m >b, \underline{S}_m >0, \tau^+_{\sqrt{m}} \leq \omega m \right)+P_a\left({S}_m >b, \underline{S}_m >0, \tau^+_{\sqrt{m}} > \omega m \right)=:P_3+P_4.
\end{align*}
Strong Markov property and homogeneity give: 
\begin{align*}
P_a\left({S}_m >b, \underline{S}_m >0, \tau^+_{\sqrt{m}}= j \right) &\leq P_a\left(\tau^-_0> \tau^+_{\sqrt{m}}= j \right) \sup_{0 \leq x \leq \alpha }P_{\sqrt{m}+x}\left({S}_{m-j} >b, \underline{S}_{m-j} >0 \right)\\
&\leq P_a\left(\tau^-_0> \tau^+_{\sqrt{m}}= j \right) P\left({S}_{m-j} >b-\sqrt{m}- \alpha \right)
\end{align*}
implying with Lemma \ref{B.1}:
\begin{align*}
 P_4 & \leq \sum_{j=\omega m+1}^{m} P_a\left(\tau^-_0> \tau^+_{\sqrt{m}}= j \right) P\left({S}_{m-j} >b-\sqrt{m}- \alpha \right) \\
& \leq P_a\left(\tau^-_0> \tau^+_{\sqrt{m}} \right) \sup_{ \omega m \leq j \leq m  } P\left({S}_{m-j} >b- \sqrt{m}-\alpha \right) \\ 
& \leq C_+ \frac{a}{\sqrt{m}} \sup_{ \omega m \leq j \leq m  } P\left({S}_{m-j} >b- \sqrt{m}-\alpha \right),
 \end{align*}
A classical result of moderate deviations (see for instance  \cite{Petrov}, Chapter VIII, Theorem 1) implies 
 \begin{align*}
 P\left({S}_{m(1- \omega)} >b- \sqrt{m}-\alpha \right) &\leq \frac{C_+  \sqrt {m (1- \omega)} }{(b-\sqrt{m}- \alpha)} e^{(b- \sqrt m- \alpha)g( (b- \sqrt m- \alpha)/(m(1- \omega)))} \\
 & \sim \frac{C_+ \sqrt m }{b} e^{(b- \sqrt m)g( (b- \sqrt m)/(m(1- \omega)))} 
 \end{align*}
We now choose $\omega$ in such a way that $(b- \sqrt m)g( (b- \sqrt m)/(m(1- \omega)))-b g(b/m)=O(1)$, $\omega$ is actually a sum which first two terms are $\omega=2 m^{1/2}/b-\lambda(b/m)/\sqrt{m}+ ...$. So for any $n$ large enough 
 \begin{align*}
 P\left({S}_{m(1- \omega)} >b- \sqrt{m}-\alpha \right) &\leq  \frac{C_+ \sqrt m }{b} e^{b \cdot g( b/m)}. 
 \end{align*}


\noindent With similar computations this upper bound is still true for $P\left({S}_{m-j} >b- \sqrt{m}-\alpha \right)$ for $m(1-\varepsilon) \leq j \leq m $, so $P_4 \leq C_+ \frac{a}{{b}} e^{b  g( b/m)} $. In the same way
\begin{align*}
 P_3 & \leq \sum_{j=\sqrt{m}/ \alpha}^{\omega m} P_a\left(\tau^-_0> \tau^+_{\sqrt{m}}= j \right) \sup_{0 \leq x \leq \alpha }P_{\sqrt{m}+x}\left({S}_{m-j} >b, \underline{S}_{m-j} >0 \right) \\
 & \leq \sum_{j=\sqrt{m}/ \alpha}^{\omega m} P\left(S_j \geq {\sqrt{m}}-a \right) P\left({S}_{m-j} >b-\sqrt{m}- \alpha \right),
 \end{align*}
 Using again  \cite{Petrov}, 
 \begin{align*}
 P_3 & \leq C_+ \sum_{j=\sqrt{m}/ \alpha}^{\omega m} \frac{j^{1/2}}{\sqrt{m}-a}e^{-m/(2 \sigma^2j )} \frac{(m-j)^{1/2}}{b-\sqrt m-\alpha}e^{(b-\sqrt m-\alpha)g((b-\sqrt m-\alpha)/(m-j) )} \\
& \leq C_+ \frac{(\omega m)^{1/2}}{b} e^{-b/(2 \sigma^2 m^{1/2} )} e^{b \cdot g( b/m)}=o\left(\frac{  e^{b \cdot g( b/m)}.}{b} \right).
 \end{align*}
 which finish the proof. 
\noindent \eqref{beautiful3} can be proved in a similar way with classical large deviation estimates.
\end{Pre} 
 
\noindent \\
The following Lemma states the local behavior of sums of i.i.d. random variables, recall that $S_1$ is non-lattice.
 
 \begin{lem} \label{beautiful} Let $\epsilon>0$ small  and $A>0$ large. For all  $m$ large enough, for all $1 \leq r \leq A  \sqrt m$
 \begin{align}
P\left({S}_m \in (r,r+1],\underline{S}_m >0 \right)= \frac{  r }{ m^{3/2}} e^{-r^2/  (2 \sigma^2 m) }+ o(m^{-3/2}). \label{B.2}
\end{align}
For all $ A  \sqrt m \leq r \leq \epsilon m$
 \begin{align}
P\left({S}_m \in (r,r+1], \underline{S}_m >0 \right)\asymp \frac{1}{m} e^{r g(r/m)}.
 \label{1.4}
\end{align}
see Lemma \ref{beautiful2} for the definition of $g$.
\end{lem}
 
\begin{Pre}
\eqref{B.2} is F. Caravenna \cite{Carav} result and  \eqref{1.4} can be obtained with \cite{Petrov} Chapter VIII, Theorem 2 and 10 and similar arguments than in the proof of Lemma \ref{beautiful2}.
 \end{Pre}
 

\section{Probability of hitting time}

\begin{Lem}\label{lemfinal}
For $x^\prime\in\llbracket \phi,x\rrbracket$:
\begin{eqnarray}
\p_{x^\prime_x}^{\mathcal E}(T_x<T_{x^\prime})&=&\frac{e^{V(x^\prime_x)}}{\sum_{z\in\rrbracket x^\prime,x\rrbracket }e^{V(z)}} \label{eq1},\\
\p_{\overset{\leftarrow}{x}}^{\mathcal E}(T_{x^\prime}<T_x)&=&\frac{e^{V(x)}}{\sum_{z\in\rrbracket x^\prime,x\rrbracket }e^{V(z)}}.\label{eq2}
\end{eqnarray}
where $x^\prime_x$ is the only children of $x^\prime$ in $\llbracket x^\prime,x\rrbracket$. 
\end{Lem}

The result is classical (see for instance \cite{AndreolettiDebs1}) and a useful direct consequence of this latter is the following
\begin{cor-anglais}\label{routine}
Let  $\mathcal A\subset \lbrace z\in\mathbb T,\vert z\vert =\ell\rbrace$ and $\kappa>0$, there exists a positive constant $c_7$ such that
\begin{align}
\pe_\phi(T_\phi^{n^{\kappa}}<T_z)&\leq \exp\left(-c_7 n^{\kappa}e^{{\overline V(z)}}/{\ell}\right),\, \forall z\in\mathcal A, \label{C3}\\
\pe\left(\bigcup_{z\in\mathcal A}\lbrace\lo(z,T_\phi^{n^{\kappa}})=0\rbrace\right)&\leq \vert \mathcal A\vert \exp\left(-c_7n^{\kappa}e^{-{\max_{z\in\mathcal A}\overline V(z)}}/{\ell}\right)\label{C4}
\end{align}
\end{cor-anglais}

\begin{Pre}
Obviously \eqref{C4} is a consequence of \eqref{C3}. 
Thanks to formula \eqref{eq1}, for $z\in\mathcal A$, $\pe_{\phi}\left(T_z<T_\phi\right) $ $\geq C_-  e^{-\overline V(z)}/\ell$. 
Then  using the strong Markov property and the recurrence of $X$, for $n$ large enough
$\pe_{\phi}\left(T_\phi^{n^\kappa}<T_z\right)=(1-\pe_{\phi}\left(T_z<T_\phi\right))^{n^{\kappa}}\leq \exp\left(e^{-c_7n^{\kappa}{\overline V(z)}}/{\ell}\right).$
\end{Pre}

\bibliographystyle{plain}
\bibliography{thbiblio}

\end{document}